\documentclass[a4paper,11pt,twoside]{article}
\usepackage{pst-fill,pst-grad}
\usepackage{textcomp}
\usepackage[english]{babel}
\usepackage[utf8]{inputenc} 
\usepackage{graphicx}
\usepackage{amsmath}
\usepackage{float}
\usepackage{fancyhdr}
\usepackage[matrix,arrow,curve]{xy}
\usepackage{pstricks} 
\usepackage{amsmath,amsfonts,verbatim,afterpage,theorem,euscript,mathrsfs,amssymb}
\usepackage{amsfonts}
\usepackage{amssymb}
\usepackage{array}
\usepackage{dsfont}
\usepackage{hyperref}
\usepackage{authblk}
\usepackage{color}
\newtheorem{Proposition}{Proposition}[section]
\newtheorem{Lemme}{Lemma}[section]
\newtheorem{Theoreme}{Theorem}[section]
\newtheorem{TheoremeP}{Theorem}

\newtheorem{Remarque}{Remark}
\newcommand{\sign}{\text{sign}}
\def \R{\mathbb{R}}

\def \finpv{\hfill $\blacksquare$  \\ \newline }
\def \pv{{\bf{Proof.}}~} 

\def \ds{\displaystyle}
\title{\bf On decay properties and asymptotic behavior of solutions to  a non-local perturbed \emph{KdV} equation}

\author[1]{ Manuel Fernando Cortez\footnote{\emph{Corresponding author}: manuel.cortez@epn.edu.ec}}
\author[2]{ Oscar Jarr\'in\footnote{or.jarrin@uta.edu.ec}}

\affil[1]{\scriptsize Departamento de Matem\'aticas, Escuela Politécnica Nacional,  Ladron de Guevera E11-253, Quito, Ecuador}
\affil[2]{\scriptsize Dirección de investigación y desarrollo (DIDE),
	Universidad Técnica de Ambato, 
	Av. de los Chasquis, 180207, Ambato, Ecuador.}
\begin{document}
\maketitle

%%%%%%%%%%%%%%%%%%%%%%%%%%%%%%%%%%%%%%%%%%%%%%
\begin{scriptsize}
\abstract{We consider the \emph{KdV} equation with an additional non-local perturbation term defined through the Hilbert transform, also known as the OST-equation. We prove that the solutions $u(t,x)$ has a pointwise decay in spatial variable:  $\ds{\vert u(t,x)\vert \lesssim \frac{1}{1 + |x|^{2}}}$, provided that the initial data has the same decaying and moreover we find the asymptotic profile of $u(t,x)$ when $|x| \to +\infty$. \\
\\
Next, we show that decay rate given above is optimal when the initial data is not a  zero-mean function and in this case we derive an estimate from below $\ds{\frac{1}{\vert x\vert^2} \lesssim \vert u(t,x)\vert}$ for $\vert x \vert$ large enough. In the case when the initial datum is a zero-mean function, we prove that the decay rate above is improved to $\ds{\frac{1}{1+\vert x \vert^{2+\varepsilon}}}$ for $0<\varepsilon \leq 1$. Finally, we study the local-well posedness of the OST-equation in the framework of  Lebesgue spaces.}\\[1cm]
\end{scriptsize}
%\tableofcontents
%%%%%%%%%%%%%%%%%%%%%%%%%%%%%%%%%%%%%%%%%%%%%%
%%%%%%%%%%%%%%%%%%%%%%%%%%%%%%%%%%%%%%%%%%%%%%
\section{Introduction} 
In this article we consider the following   Cauchy's problem for a non-locally perturbed \emph{KdV} equation
\begin{equation}  \label{eq:f}
\left\{
\begin{array} [c]{l} \vspace{2mm}  %
\partial_t u     + u \partial_x u + \partial^{3}_{x} u  + \eta(\mathcal{H} \partial_x u  + \mathcal{H} \partial^{3}_{x} u ) = 0, \quad \eta>0, \quad \text{on}\quad ]0,+\infty[\times \mathbb{R},\\  
u(0,\cdot) = u_0.
\end{array}
\right. 
\end{equation}
where the function $u:[0,+\infty[\times \mathbb{R}   \rightarrow \mathbb{R}$ is the solution, $u_0: \mathbb{R} \rightarrow \mathbb{R}$ is the initial datum and $\mathcal{H}$ is the Hilbert transform defined as follows: for $\varphi \in \mathcal{S}(\mathbb{R})$,
\begin{equation}
\label{HIL}
\mathcal{H}(\varphi)(x)=p.v. \frac{1}{\pi} \int_{\mathbb{R}} \frac{\varphi(x-y)}{y} dy.
\end{equation} 
Equation (\ref{eq:f}), also called the Ostrovsky, Stepanyams and Tsimring equation (OST-equation), was derived by Ostrovsky \emph{et al.} in 
 \cite{OST,OST1} to describe the  radiational  instability of long non-linear waves in a stratified flow caused by internal wave radiation from a shear layer. It  deserves  remark that when $\eta=0$ we obtain the well-know \emph{KdV} equation. The parameter $\eta>0$ represents the importance of amplification and damping relative to dispersion. Indeed, the fourth term in equation (\ref{eq:f}) represents amplification, which is responsible for the radiational instability of the negative energy wave, while the fifth term in equation (\ref{eq:f}) denotes damping (see \cite{OST0} for more details).  Both of these two terms are described by the non-local integrals represented by the Hilbert transform \eqref{HIL}.\\
\\
One rewrites Equation \eqref{eq:f} in the equivalent Duhamel formulation (see \cite{BorysAlvarez-tesis}):
\begin{equation}\label{integral}
u(t,x)=K_\eta(t,\cdot)\ast u_0(x)-\frac{1}{2} \int_{0}^{t}K_{\eta}(t-\tau,\cdot)\ast \partial_{x}(u^{2})(\tau,\cdot)(x)d \tau, 
\end{equation}  where the kernel $K_\eta(t,x)$  is given by 
\begin{equation}\label{Kernel} 
K_{\eta}(t,x)= \mathcal{F}^{-1}\left( e^{(i\xi^3t - \eta t(\vert \xi \vert^3-\vert \xi \vert))} \right)(x), 
\end{equation} and where $\mathcal{F}^{-1}$ denotes the inverse Fourier transform.\\
\\
Well-posedness results for the Cauchy problem (\ref{integral}) was extensively studied in the framework of Sobolev spaces. The first work on this problem was carried out by B. Alvarez Samaniego in his PhD thesis \cite{BorysAlvarez-tesis} (see also the article \cite{BorysAlvarez-2} by the same author). Alvarez Samaniego  proved the local well-posedness  in $H^{s} (\mathbb{R})$ for $s>\frac{1}{2}$, using properties of the semi-group associated with the linear problem. He also obtained a global solution in $H^{s}$ for $s\geq 1$, making use of the standard energy estimates.
This result was improved by several authors: X. Carvajal \& M. Scialon  proved in their article  \cite{XCarvajal-MScialom}, through Strichartz-type  estimates and smoothing effects, the local well-posedness (LWP) of the Cauchy's problem (\ref{integral}) in $H^{s}(\R)$ for $s\geq 0$ and global well-posedness (GWP) in $L^{2}(\R)$.  After, X. Zhao \& S. Cui proved in \cite{ZhaoCui} and \cite{ZhaoCui-2} the LWP  of problem (\ref{eq:f}) in $H^{s}(\R)$ for $s>-\frac{3}{4}$ and GWP for $s\geq 0$. Finally, in recent works A. Esfahani and H. Wang \cite{AEsfahani, AEsfahani1} used purely  dissipative approaches based on the  method of  bilinear estimates in the Bourgain-type spaces (see also \cite{Molinet} for more references on these spaces) to show  that the Cauchy problem  (\ref{integral}) is LWP in $H^{s}(\R)$ for $s \geq -\frac{3}{2}$ and moreover, it is shown that  $H^{-\frac{3}{2}}$ is the critical Sobolev space for the LWP.\\
\\
On the other hand, since equation (\ref{eq:f}) is a nonlinear  dissipative  equation, it is natural to ask for the existence of solitary waves. \emph{Numerical} studies 
done in \cite{Feng} by B.F.  Feng  \&  T.  Kawahara shows that  for   every $ \eta>0$ there exists  a  family  of   solitary  waves  which \emph{experimentally} decay  as $\ds{\frac{1}{1+ |x|^{2}}}$ when $|x| \to + \infty$.  This numerical decay of solitary waves suggests the theoretical study of the decay in spacial variable of  solutions $u(t,x)$ of equation (\ref{eq:f}) and, in this setting,  B. Alvarez Samaniego showed  in the last part of his PhD thesis (see Theorem $5.2$ of \cite{BorysAlvarez-tesis})  that if the initial datum $u_0$ verifies $\ds{u_0 \in H^{2}(\R) \cap L^{2}(1+ \vert \cdot \vert^{2}, dx)}$ then there exists  $\ds{u \in \mathcal{C}([0, \infty[,H^{2}(\R) \cap L^{2}(1+ \vert \cdot \vert^{2}, dx))}$ a unique solution of equation (\ref{eq:f}). This result is intrinsically related to the nature of the functional spaces above  in which the Fourier Transform plays a very important role: kernel $K_{\eta}(t,x)$ given in (\ref{Kernel})  associated with the equation is explicitly defined in frequency variable. Furthermore, remark that this spatially-decaying of solution  is studied in  the setting  of the  weighted-$L^{2}$ space and therefore it's a weighted average  decay. \\
\\
The general aim of this paper is to study  spatial decay estimates of the solution $u(t,x)$. Our methods are inspired by L. Brandolese \emph{et. al.}  \cite{Bra,Bra2} which  are  essentially based on well-know properties of the kernel associated to the linear equation, however, our approach to find these estimates  is a little different.  Indeed, using the explicit definition in the frequency variable of $K_\eta(t,x)$ and the inverse Fourier transform,  we deduce some sharp spatially-decaying properties  for this kernel and for its derivatives.  It is worth  remarking that these methods are technically different with respect to previous works on equation \eqref{eq:f} since in those works the kernel is studied in the frequency variable and not in the spatial variable. \\
\\
On the other hand,  this approach permits to study the equation \eqref{eq:f} in other functional spaces which, to the best of our knowledge, have not been considered before. More precisely, we prove that the  properties in the spatial variable of  kernel $K_\eta(t,x)$ allow us to prove that the integral equation (\ref{integral}) is LWP  for small initial datum in the framework of Lebesgue spaces. \\
\\
\textbf{Organization of the paper.}  In Section \ref{sec:Results} we state all the results obtained. In Section \ref{sec:kernel} we study the optimal decay in spatial variable of the kernel $K_{\eta}(t,x)$.  Section \ref{sec:pointwise-decay-asymptotics} is devoted to the study of pointwise decaying and asymptotic behavior in the spacial variable of solutions of equation (\ref{eq:f}).  The last section \ref{sec:Lebesgue} is devoted to the studies of the LWP of equation (\ref{eq:f}) in the framework of Lebesgue spaces.\\
\\
\textbf{Acknowledgments.} We thank the referees  for the useful remarks and  comments which allow us to improve our work.
\section{Statement of the results}\label{sec:Results} 
\subsection{Pointwise decay and asymptotic behavior in spatial variable} 
The first purpose of this paper is to obtain a \emph{pointwise}  decay in the spacial variable of solution $u(t,x)$. More  precisely, we prove that if the initial datum $u_0 \in H^{s}(\mathbb{R})$  (with $\frac{3}{2}<s\leq 2$)  verifies $\ds{\vert u_0(x)\vert \leq \frac{c}{1+\vert x \vert^2}}$, then there exist a unique global in time solution  $u(t,x)$ of the integral  equation   (\ref{integral}) which  fulfills  the same decay of the  initial datum $u_0$. Moreover, we show that the solution $u(t,x)$ of the integral equation (\ref{integral})  is smooth enough  and then this solution verifies the differential equation (\ref{eq:f}) in the classical sense.
\begin{TheoremeP}\label{Th-decay}   Let $\frac{3}{2} < s\leq 2$ and  let $u_0 \in H^{s} (\R)$ be an initial datum, such that  $\ds{\vert u_0(x)\vert \leq \frac{c}{1+\vert x \vert^{2}}}$. Then, the equation  \eqref{eq:f} possesses a unique  solution $u \in \mathcal{C}(]0,+\infty[, \mathcal{C}^{\infty}(\R))$ arising from $u_0$, such that for all time $t>0$ there exists a constant $C(t,\eta,u_0,u)>0$ such that for all $x\in \R$ the solution $u(t,x)$ verifies:  
\begin{equation}\label{decay}	
\vert u(t,x)\vert \leq \frac{C(t,\eta,u_0,\|u\|_{H^{s}})}{1+\vert x \vert^2}. 
\end{equation}
\end{TheoremeP}	 
\begin{Remarque} Estimate  $(\ref{decay})$ is valid only in the setting of the perturbed \emph{KdV} equation (\ref{eq:f})  when the parameter $\eta$ is strictly positive.
\end{Remarque}	
Indeed, with respect to the parameter $\eta$ the constant $C(t,\eta,u_0,u)>0$ behaves like the following  expression (see formula (\ref{Constante-C}) for all the details):  $\ds{\frac{1}{\eta^{\frac{1}{3}}} \left( 1 + \left(\frac{1}{\eta} +2 \right)^{2} \right) +1}$,  and this expression  is not controlled when  $\eta \longrightarrow 0^{+}$.\\
\\
Recall that in the case $\eta=0$ the equation  \ref{eq:f} becomes the \emph{KdV} equation.  In this framework  T. Kato \cite{Kato} showed the following  persistence problem:   if $u_0 \in H^{2m}(\R) \cap L^{2}(|x|^{2m}, \,dx)$, where   $m\in \mathbb{N}$ is strictly positive, then  the Cauchy problem for the \emph{KdV} equation is globally well-posed in the space $C([0,+ \infty[;  H^{2m}(\R) \cap L^{2}(|x|^{2m} \,dx))$ and then the solution of the \emph{KdV} equations decays at infinity as fast as the initial datum. For related results see also \cite{FonsecaLinaresPonce} and \cite{NahasPonce}.  \\
\\
Getting back to the perturbed \emph{KdV} equation (\ref{eq:f})  a natural question arises: is the spatial decay given in the formula (\ref{decay}) optimal? and concerning this question  B Alvarez Samaniego has shown in \cite{BorysAlvarez-1} that the solution cannot have a weight average decay faster than $\ds{\frac{1}{1 + | x |^{4}}}$; and in this case we have a loss of persistence in the spatial decay. This results suggests  that the optimal  decay rate in spatial variable of solution $u(t,x)$ must be of the order $\ds{\frac{1}{1+ | x|^{s}}}$ with $2 \leq s < 4$. \\
\\
The second purpose of this paper is to study how sharp is the decay rate of solution given in Theorem \ref{Th-decay}. For this purpose,  in the following theorem we start by studying the asymptotic profile of solution $u(t,x)$ and we prove that  if the initial datum  $u_0$ decays a little faster than $\ds{\frac{1}{ 1+| x|^{2}}}$, then the solution $u(t,x)$ associated to $u_0$ has the following asymptotic behavior in the spatial variable.   
\begin{TheoremeP}\label{Th:asymptotics} Let $\frac{3}{2} < s\leq 2$ and  let $u_0 \in H^{s} (\R)$ be an initial datum such that for $\varepsilon>0$ we have $\ds{\vert u_0(x)\vert \leq \frac{c}{1+\vert x \vert^{2+\varepsilon}}}$ and  $\ds{ \left\vert \frac{d}{dx} u_0(x)\right\vert \leq \frac{c}{1+\vert x \vert^2}}$. Then, the solution $u(t,x)$ of the equation (\ref{eq:f})  given  by Theorem \ref{Th-decay}  has the following asymptotic development when $\vert x \vert$ is large enough:
\begin{equation}\label{asymptotic}
u(t,x)= K_\eta(t,x) \left( \int_{\R}u_0(y)dy\right)+ \int_{0}^{t} K_\eta(t-\tau, x) \left( \int_{\R}u(\tau,y) \partial_y u(\tau,y) dy \right) d \tau +\, o(t)\left(  1 / \vert x \vert^2\right), 
\end{equation} where the kernel $K_{\eta}(t,x)$ is given in (\ref{Kernel}), and where the quantity  $o(t)\left( 1/ \vert x \vert^2 \right)$ is such that for all $t>0$
\begin{equation}\label{little-o} 
 \lim_{\vert x \vert \longrightarrow +\infty}  \frac{o(t)(1/\vert x \vert^2)}{1/\vert x \vert^2 }=0.
 \end{equation}
 \end{TheoremeP}
This  asymptotic development of solution $u(t,x)$ provides us interesting information on the behavior of this solution in spatial variable. Remark first that  all the information respect to the spatial variable relies on the information (in the spatial variable) of the kernel $K_\eta(t,x)$. More precisely, concentrating our attention in the first term on the right-hand side of this identity we may observe that this term is not zero when the initial datum $u_0$ is not a zero-mean function ($\int_{\mathbb{R}}u_0(y)dy \neq 0$). Moreover,  in Proposition \ref{Prop1} below, we show that this kernel has an optimal decay rate of the order $\frac{1}{1+\vert x \vert^2}$ and this fact suggests that the decay of solution $u(t,x)$ given in Theorem \ref{Th-decay} must be sharp when the initial datum verifies a non zero-mean condition. \\ 
\\ 
On the other hand, in the case of a zero-mean initial datum ($\int_{\mathbb{R}} u_0(y) dy=0$) and for $\vert x \vert $ large enough, observe that the solution behaves essentially  as the second term on the right-hand side of  identity  \eqref{asymptotic} which comes from  the nonlinear term in equation (\ref{integral}).  In this case we shall prove that the decay rate of solution given in the formula (\ref{decay}) actually  is not sharp and it can be improved. \\
\\
Our next result summarizes these statements. 
\begin{TheoremeP}\label{Th:estim-below}  Under the same hypothesis of   Theorem  \ref{Th:asymptotics}. 
\begin{enumerate}
\item[1)]  Assume that  $\ds{\int_{\mathbb{R}}u_0(y)dy \neq 0}$. Then there exists  $M>0$ and there exists a constant $0<c_{\eta, t} <c_\eta e^{4 \eta t}$ such that for $\vert x \vert >M$ we have the estimate from below: 
\begin{equation}\label{estim-below} 
\frac{c_{\eta,t} \, t}{2 \vert x \vert^2}  \left\vert \int_{\mathbb{R}}u_0(y)dy \right\vert \leq   \vert u (t,x) \vert. 
\end{equation}
\item[2)]  Assume that $\ds{\int_{\mathbb{R}}u_0(y)dy = 0}$. Then the solution $u(t,x)$ of the equation (\ref{eq:f})  given  by Theorem \ref{Th-decay} has the following decay: for $0 < \varepsilon \leq 1$ 
\begin{equation}\label{decay2}
\vert u(t,x)\vert \leq \frac{C'(\eta,\varepsilon,t,u_0,u)}{1+\vert x \vert^{2+\varepsilon}}, 
\end{equation} where the constant  $C'(\eta,\varepsilon, t, u_0, \|u\|_{H^{s}})>0$ does not depend on the variable $x$. 
\end{enumerate}		
\end{TheoremeP}	
%REVISADO GRAMATICA
\begin{Remarque}
It should be emphasized that while the non zero-mean condition $\int_{\mathbb{R}}u_0(y)dy \neq 0$ is verified, even if the initial datum is a smooth, compact-support  function  the arising solution $u(t,x)$ cannot decay at  infinity faster than $\ds{\frac{1}{|x|^2}}$ and in this case the decay rate of solution given in Theorem \ref{Th-decay} is optimal.
\end{Remarque}
% \\
%\\
\begin{Remarque}
When $\displaystyle\int_{\mathbb{R}}u_0(y)dy = 0$ , the estimate from below (\ref{estim-below}) is  not more valid and moreover the  decay rate of solution $u(t,x)$ given in Theorem \ref{Th-decay} is  improved in estimate (\ref{decay2}). Thus,  the persistence problem is valid for $0< \epsilon \leq 1$.  
\end{Remarque}
However, with  respect to optimality of this estimate, our approaches  do not seem to be sufficient to derive an estimate from below of the type $\frac{1}{\vert x \vert^{2 + \epsilon}} \lesssim \vert u (t,x)\vert$. This fact remains an interesting open question. 
\begin{Remarque}
For $\epsilon>1$, the persistence problem studied in point $2)$ of Theorem \ref{Th:estim-below} does not seem to be valid.
\end{Remarque}
Indeed, roughly speaking,  inequality (\ref{decay2})  relies on sharp estimates for the linear and the nonlinear term  in the integral formulation of the solution given in (\ref{integral}). The estimate done on the linear term actually can be improved as
$$ \vert K_\eta (t,x)\ast u_0(x) \vert \lesssim \frac{1}{1+ |x|^{2 + \epsilon}},$$ with $\epsilon>1$,  provided that the initial datum is a zero-mean function which decays fast enough, but, the nonlinear term is estimated as
$$ \left\vert   \int_{0}^{t}\partial_{x} K_{\eta}(t-\tau,\cdot)\ast u^{2}(\tau,\cdot)(x)d \tau \right\vert \lesssim \frac{1}{1+\vert x \vert^3},$$ see estimate (\ref{nonlinterm}) for the details.  As  the expression  $\partial_x K_\eta(t,x)$ has a sharp decay of the order $\frac{1}{\vert x \vert^3}$  this  proposes that this term cannot decay faster than $\frac{1}{\vert x \vert^3}$ and, to the best of our knowledge, we do not know a better estimate. 

\subsection{The local well-posedness in Lebesgue spaces} 

The third purpose of this paper is  to study the existence and uniqueness of \emph{mild} solutions for the Cauchy problem (\ref{eq:f}) in the framework of Lebesgue spaces when the initial datum $u_0$ is small enough. We start by recalling that we refer to a \emph{mild} solution $u(t,x)$ when this solution is written as the integral formulation (\ref{integral}) and it is obtained by a fixed-point argument.\\
\\
It is worth  remarking here that the following theorem is just a first study in the setting of Lebesgue spaces and we think that this result could be improved in further investigations. 
\begin{TheoremeP}\label{Th-Lebesgue} Let $1\leq p \leq +\infty$ and let $u_0 \in L^{p}(\R)$ be an initial datum. Let $T>0$. Then, there exists $\delta=\delta(T)>0$ such that if $\Vert u_0 \Vert_{L^p}<\delta$ then the integral equation (\ref{integral}) possesses a unique  solution local in time  $u \in L^{\infty}(]0, T[, L^{p}(\R))$ which verifies  $\ds{\sup_{0\le t \leq T} t^{\frac{1}{3}}\Vert u(t,\cdot)\Vert_{L^p}<+\infty }$.  \\
\end{TheoremeP}
\begin{Remarque}
The value of the parameter $p=2$ is of particular interest since in this case the result above gives a new proof for the LWP obtained by X. Carvajal $\&$ M. Scialons in \cite{XCarvajal-MScialom}, which relies essentially on smoothing effects and Strichartz-type  estimates.
\end{Remarque}
We finish the statement of our results with the following interesting remark.
\begin{Remarque}
All our results given for the equation \eqref{eq:f} are still valid (under some technical modifications in the proofs) for the non-local perturbation of the Benjamin-Ono (npBO) equation:
\begin{equation}  \label{npBO}
\left\{
\begin{array} [c]{l} \vspace{2mm}  %
\partial_t u     + u \partial_x u + \mathcal{H} \partial^{2}_{x} u  + \eta(\mathcal{H} \partial_x u  + \mathcal{H} \partial^{3}_{x} u ) = 0, \quad \eta>0, \quad \text{on}\quad ]0,+\infty[\times \mathbb{R},\\  
u(0,\cdot) = u_0.
\end{array}
\right. 
\end{equation}
\end{Remarque}
 Note first that the only difference between equation \eqref{eq:f} and \eqref{npBO} is the  linear term $\ds{\mathcal{H} \partial^{2}_{x} u }$. This equation was recently studied by Foseca \emph{et. al.} in \cite{FonPas}, where they proved similar results concerning the local and global well-posedness, and regularity issues.\\
\\
On the other hand, it is easy to see that the  kernel $F_\eta(t,x)$ associated to  equation \eqref{npBO} in explicitly defined in the  frequency variable as
\begin{equation*}
\widehat{F_\eta} (t, \xi)= e^{i \xi |\xi|  t + t \eta(|\xi| - |\xi|^{3})},
\end{equation*}
and thus, our methods can be adapted without any problem. Indeed, our results are purely based on estimates on the non-complex exponential part $e^{ t \eta(|\xi| - |\xi|^{3})}$ which is exactly the same for the kernel  $\widehat{K_\eta} (t, \xi)$.
%%%%%%%%%%%%%%%%%%%%%%%%%%%%%%%%%%%%%%%%%%%%%%%%%%%%%%%%%%%%%%%%%%%%%%%%%%%%%%%%%%%%%%%%%%
\section{kernel estimates}\label{sec:kernel}
In this section we study the properties decay in spacial variable of the kernel $K_{\eta}(t,x)$ which will be useful in the next sections.
\begin{Proposition}\label{Prop1}
	Let $K_\eta(t,x)$ be the kernel defined in the expression (\ref{Kernel}). 
	\begin{enumerate}
		\item[1) ]  There exists a constant $c_\eta>0$,  given in the  formula  (\ref{Const-eta}) and which only depends on $\eta>0$, such that for all time $t>0$ we have $\ds{\vert K_\eta(t,x)\vert \leq  c_\eta \frac{e^{5\eta t}}{t^{\frac{1}{3}}} \frac{1}{1+\vert x \vert^2}}$.
		\item[2)] Moreover,  the kernel $K_{\eta}(t,x)$ cannot decay at infinity faster than $\ds{\frac{1}{1+\vert x \vert^2}}$.  
	\end{enumerate}	
	%Then, there exists a constant $c_\eta>0$, which only depends of $\eta>0$, such that for all time $T>0$ we have $$  \sup_{t \in [0,T]} t^{\frac{1}{3}} \Vert (1+\vert \cdot \vert^2) K_\eta(t,\cdot)\Vert_{L^{\infty}} \leq c_\eta e^{5\eta T}.$$ 
\end{Proposition} 
\pv
\begin{enumerate}
	%revisado parte del kernel 23-033
	\item[1)] First, we will estimate the quantity $\vert K_\eta (t,x)\vert$ and then we will estimate the quantity $\vert x \vert^{2} \vert K_{\eta}(t,x)\vert.$\\
	\\
	We write 
	\begin{equation}\label{eq09}
	\vert K_{\eta}(t,x)\vert \leq \Vert K_{\eta}(t,\cdot)\Vert_{L^{\infty}}\leq \Vert \widehat{K_{\eta}}(t,\cdot)\Vert_{L^1},
	\end{equation}  and then we must study  the term $\Vert \widehat{K_{\eta}}(t,\cdot)\Vert_{L^1}$.    By the expression (\ref{Kernel}), we have $\widehat{K_{\eta}}(t,\xi)= e^{(i\xi^3 t - \eta t(\vert \xi \vert^3-\vert \xi \vert))}$ and  we can write 
	\begin{eqnarray}\label{eq02} \nonumber
	\Vert \widehat{K_{\eta}}(t,\cdot)\Vert_{L^1} &=& \int_{\mathbb{R}}  \vert e^{i\xi^3 t} \vert \vert  e^{ - \eta t(\vert \xi \vert^3-\vert \xi \vert)} \vert d \xi  = \int_{\mathbb{R}}  e^{ - \eta t(\vert \xi \vert^3-\vert \xi \vert)} d \xi \\ \nonumber
	& =&  \int_{\vert \xi \vert \leq \sqrt{2}}  e^{ - \eta t(\vert \xi \vert^3-\vert \xi \vert)} d \xi 
	+ \int_{\vert \xi \vert > \sqrt{2}}  e^{ - \eta t(\vert \xi \vert^3-\vert \xi \vert)} d \xi  \\ 
	&=& I_1+I_2.  
	\end{eqnarray} 
	In order to estimate the integral $I_1$, remark that if $\vert \xi \vert \leq \sqrt{2}$ then we have $-(\vert \xi \vert^3 -\vert \xi \vert) \leq \vert \xi \vert$ and thus we can write 
	$$ I_1 \leq \int_{\vert \xi \vert \leq \sqrt{2}} e^{\eta t \vert \xi \vert} d \xi \leq c\, e^{ \sqrt{2}\eta t} \leq c\,e^{2 \eta t}.$$ 
	Now, in order to estimate the integral $I_2$, remark that if $\vert \xi \vert >\sqrt{2}$ then we have $-(\vert \xi \vert^3-\vert \xi \vert)<-\frac{\vert \xi \vert^3}{2}$ and thus, we write 
	$$ I_2 \leq \int_{\vert \xi \vert > \sqrt{2}} e^{-\eta t \frac{\vert \xi \vert^{3}}{2}} d \xi \leq \int_{0}^{+\infty}  e^{-\eta t \frac{\vert \xi \vert^{3}}{2}} d \xi  \leq \frac{c}{(\eta t)^{\frac{1}{3}}}.$$
	%finpagina6
	With these estimates, we get back to  the identity  (\ref{eq02}) and we write  
	\begin{equation}\label{Estim-K-Fou-L1}
	\Vert \widehat{K_{\eta}}(t,\cdot) \Vert_{L^1} \leq c\,e^{2 \eta t} +\frac{c}{(\eta t)^{\frac{1}{3}}}\leq c\, \frac{e^{2 \eta t} (\eta t)^{\frac{1}{3}} +1}{(\eta t)^{\frac{1}{3}}} \leq c\, \frac{e^{3\eta t}+1}{(\eta t)^{\frac{1}{3}}}  \leq C \frac{e^{3\eta t }}{(\eta t)^{\frac{1}{3}}},
	\end{equation} hence,   getting back to  the estimate (\ref{eq09}), we can write 	
	\begin{equation}\label{eq03}
	\vert K_\eta(t,x)\vert \leq  C \frac{e^{3\eta t }}{(\eta t)^{\frac{1}{3}}}.  
	\end{equation}
	Now we will estimate  the quantity  $\vert x \vert^{2}  \vert K_\eta(t,x)\vert$. Recalling the expression (\ref{Kernel}), for $x\neq 0$ we write 
	\begin{eqnarray} \label{eq04} \nonumber 
	K_\eta(t,x) &=&  \mathcal{F}^{-1}\left(  e^{(i\xi^3 - \eta t(\vert \xi \vert^3-\vert \xi \vert))}  \right)(x) = \int_{\mathbb{R}} e^{2\pi i x \xi}  e^{(i\xi^3 - \eta t(\vert \xi \vert^3-\vert \xi \vert))}  d \xi \\ \nonumber  
	&=& \int_{\xi <0} e^{2\pi i x \xi}  e^{(i\xi^3 - \eta t(\vert \xi \vert^3-\vert \xi \vert))}  d \xi + \int_{\xi >0} e^{2\pi i x \xi}  e^{(i\xi^3 - \eta t(\vert \xi \vert^3-\vert \xi \vert))}  d \xi\\ \nonumber
	&=& \int_{\xi <0} e^{2 \pi i x \xi }e^{i t \xi^3- \eta t (-\xi^3 +\xi)} d \xi + \int_{\xi >0} e^{2 \pi i x \xi }e^{i t \xi^3- \eta t (\xi^3 -\xi)} d \xi \\
	&=& \frac{1}{2 \pi i x } \int_{\xi <0} 2\pi i x e^{2 \pi i x \xi }e^{i t \xi^3- \eta t (-\xi^3 +\xi)} d \xi  + \frac{1}{2\pi i x}  \int_{\xi >0} 2\pi i x e^{2 \pi i x \xi }e^{i t \xi^3- \eta t (\xi^3 -\xi)} d \xi. 
	\end{eqnarray}	
	%\end{document}
	In the last identity, remark that $\partial_{\xi}(e^{2\pi i x \xi})= 2\pi i x e^{2\pi i x \xi}$ and then, we can write 
	\begin{eqnarray*}
		& & \frac{1}{2 \pi i x } \int_{\xi <0} 2\pi i x e^{2 \pi i x \xi }e^{i t \xi^3- \eta t (-\xi^3 +\xi)} d \xi  + \frac{1}{2\pi i x}  \int_{\xi >0} 2\pi i x e^{2 \pi i x \xi }e^{i t \xi^3- \eta t (\xi^3 -\xi)} d \xi \\
		&=& \frac{1}{2 \pi i x } \int_{\xi <0} \partial_{\xi}(e^{2\pi i x \xi}) e^{i t \xi^3- \eta t (-\xi^3 +\xi)} d \xi + \frac{1}{2\pi i x}  \int_{\xi >0} \partial_{\xi}(e^{2\pi i x \xi})e^{i t \xi^3- \eta t (\xi^3 -\xi)} d \xi.
	\end{eqnarray*}
	Now, integrating by parts, each term above and  since $\ds{\lim_{\xi \longrightarrow - \infty} e^{i t \xi^3- \eta t (-\xi^3 +\xi)}=0}$ and $\ds{\lim_{\xi \longrightarrow + \infty} e^{i t \xi^3- \eta t (\xi^3 -\xi)}=0}$ then, we have 
	\begin{eqnarray*}
		& &\frac{1}{2 \pi i x } \int_{\xi <0} \partial_{\xi}(e^{2\pi i x \xi}) e^{i t \xi^3- \eta t (-\xi^3 +\xi)} d \xi + \frac{1}{2\pi i x}  \int_{\xi >0} \partial_{\xi}(e^{2\pi i x \xi})e^{i t \xi^3- \eta t (\xi^3 -\xi)} d \xi \\  \nonumber
		&=& \frac{1}{2\pi i x } - \frac{1}{2 \pi i x} \int_{\xi <0} e^{2\pi i x \xi} \partial_{\xi} \left( e^{i t \xi^3- \eta t (-\xi^3 +\xi)} \right) d \xi - \frac{1}{2\pi i x} - \frac{1}{2\pi i x}  \int_{\xi >0} e^{2\pi i x \xi} \partial_{\xi}\left( e^{i t \xi^3- \eta t (\xi^3 -\xi)}\right) d \xi \\ \nonumber
		&=& -\frac{1}{2 \pi i x} \int_{\xi <0} e^{2\pi i x \xi} \partial_{\xi} \left( e^{i t \xi^3- \eta t (-\xi^3 +\xi)} \right) d \xi -   \frac{1}{2\pi i x}  \int_{\xi >0} e^{2\pi i x \xi} \partial_{\xi}\left( e^{i t \xi^3- \eta t (\xi^3 -\xi)}\right) d \xi=(a).	
	\end{eqnarray*}  
	Thus, following the same computation done in identity (\ref{eq04}) and since   $\partial_{\xi}(e^{2\pi i x \xi})= 2\pi i x e^{2\pi i x \xi}$,  then we write  
	\begin{eqnarray}\label{eq05} \nonumber
	\ \nonumber
	(a)&=&  -\frac{1}{(2 \pi i x)^2} \int_{\xi <0} \partial_{\xi} (e^{2\pi i x \xi})  \partial_{\xi} \left( e^{i t \xi^3- \eta t (-\xi^3 +\xi)} \right) d \xi -   \frac{1}{(2\pi i x)^2}  \int_{\xi >0} \partial_{\xi}(e^{2\pi i x \xi} )\partial_{\xi}\left( e^{i t \xi^3- \eta t (\xi^3 -\xi)}\right) d \xi\\ \nonumber
	&=& -\frac{1}{(2 \pi i x)^2} \int_{\xi <0} \partial_{\xi} (e^{2\pi i x \xi}) (e^{i t \xi^3- \eta t (-\xi^3 +\xi)})(3i t \xi^2 -\eta t (-3\xi^2+1))  d \xi  \\ \nonumber 
	& & -  \frac{1}{(2\pi i x)^2}  \int_{\xi >0} \partial_{\xi}(e^{2\pi i x \xi} ) (e^{i t \xi^3- \eta t (-\xi^3 +\xi)}) (3 i t \xi^2 -\eta t (3\xi^2 -1))d \xi. \\
	&=& I_1 + I_2,
	\end{eqnarray} where we will estimate  both expressions $I_1$ and $I_2$.  For expression $I_1$, remark that we have $$ \ds{\lim_{\xi \longrightarrow - \infty} (e^{i t \xi^3- \eta t (-\xi^3 +\xi)})(3i t \xi^2 -\eta t (-3\xi^2+1))=0},$$ and integrating by parts, we can write  
	\begin{eqnarray}
	\label{eq06} \nonumber
	I_1&=&-\frac{1}{(2 \pi i x)^2}  \left( - \eta t - \int_{\xi <0} e^{2\pi i x \xi} \partial_{\xi} \left( (e^{i t \xi^3- \eta t (-\xi^3 +\xi)})(3i t \xi^2 -\eta t (-3\xi^2+1)) \right) d \xi \right)\\ %\nonumber
	&=& \frac{\eta t }{(2\pi i x)^2} + \frac{1 }{(2\pi i x)^2} \underbrace{\int_{\xi <0} e^{2\pi i x \xi} \partial_{\xi} \left( (e^{i t \xi^3- \eta t (-\xi^3 +\xi)})(3i t \xi^2 -\eta t (-3\xi^2+1)) \right) d \xi}_{= \,I_a}.  
	\end{eqnarray} Now, for the expression $I_2$ given in (\ref{eq05}), remark that we have 
	$$ \lim_{\xi \longrightarrow + \infty} (e^{i t \xi^3- \eta t (-\xi^3 +\xi)}) (3 i t \xi^2 -\eta t (3\xi^2 -1))=0,$$
	%%finpagina7
	and then, always by integration by parts  we write 
	%\end{document}
	\begin{eqnarray}\label{eq07} \nonumber 
	I_2&=& - \frac{1}{(2 \pi i x)^2}\left( - \eta t -  \int_{\xi >0} e^{2\pi i x \xi} \partial_{\xi} \left( (e^{i t \xi^3- \eta t (-\xi^3 +\xi)}) (3 i t \xi^2 -\eta t (3\xi^2 -1)) \right) d \xi \right)\\ 
	&=& \frac{\eta t }{(2\pi i x)^2} + \frac{1}{(2\pi i x)^2} \underbrace{\int_{\xi >0} e^{2\pi i x \xi} \partial_{\xi} \left( (e^{i t \xi^3- \eta t (-\xi^3 +\xi)}) (3 i t \xi^2 -\eta t (3\xi^2 -1)) \right) d \xi.}_{=\, I_b} 
	\end{eqnarray}
	Thus, with  identities (\ref{eq06}) and (\ref{eq07}) at hand, we get back to the identity (\ref{eq05}) and we write 
	\begin{equation}\label{Descomposition-Kernel}
	I_1+I_2 = \frac{2 \eta t}{(2 \pi i x)^2} +\frac{1}{(2 \pi i x)^2} (I_a + I_b),
	\end{equation} and then, getting back to the identity (\ref{eq04}) we have 
	\begin{equation}\label{eq08}
	\vert K_{\eta}(t,x) \vert =  \left\vert \frac{2 \eta t}{(2 \pi i x)^2} +\frac{2}{(2 \pi i x)^2} (I_a + I_b) \right\vert \leq c \frac{ \eta t}{x^2} + \frac{c}{x^2}  \vert I_a+ I_b \vert.  
	\end{equation} 
	We still need to estimate the term $\vert I_a+I_b\vert$ above and for this we have the following technical lemma, which we will in prove later in the appendix.
	\begin{Lemme}\label{Lemma-tech1}  There exists a numerical constant $c>0$, which does not  depend on $\eta>0$, such that for all $t>0$, we have $\ds{\vert I_a+I_b\vert \leq c \left(\frac{1}{\eta} +2 \right)^{2} e^{4 \eta t}}$.
	\end{Lemme}	
	With this estimate, we get back to the equation (\ref{eq08}) and we get 
	\begin{eqnarray*}
		\vert K_\eta(t,x)\vert & \leq & c\,\frac{\eta t}{x^{2}}+ c \left(\frac{1}{\eta} +2 \right)^{2} \frac{e^{4 \eta t} }{x^2}  \leq c \left(\frac{1}{\eta} +2 \right)^{2} \frac{\eta t }{x^2} + c \left(\frac{1}{\eta} +2 \right)^{2} \frac{e^{4 \eta t} }{x^2}\\
		& \leq & c \left(\frac{1}{\eta} +2 \right)^{2} \frac{e^{4 \eta t}}{x^2} + c \left(\frac{1}{\eta} +2 \right)^{2} \frac{e^{4 \eta t}}{x^2} \leq C \left(\frac{1}{\eta} +2 \right)^{2} \frac{e^{4 \eta t}}{x^2}.
	\end{eqnarray*} Hence, we can write 
	\begin{equation}\label{eq10}
	\vert x \vert^2 \vert K_\eta(t,x)\vert \leq C \left(\frac{1}{\eta} +2 \right)^{2} e^{4 \eta t}.
	\end{equation}
	Thus, with estimates (\ref{eq03}) and (\ref{eq10}), we can write 
	\begin{eqnarray*}
		\vert K_\eta(t,x)\vert + \vert x \vert^{2}\vert K_\eta(t,x)\vert & \leq & C \frac{e^{3\eta t }}{(\eta t)^{\frac{1}{3}}} +  C \left(\frac{1}{\eta} +2 \right)^{2} \frac{e^{4 \eta t}}{x^2} \
		\leq    C \frac{e^{3\eta t }}{(\eta t)^{\frac{1}{3}}} +  C \left(\frac{1}{\eta} +2 \right)^{2} (\eta t)^{\frac{1}{3}}\frac{e^{4 \eta t}}{(\eta t)^{\frac{1}{3}}}\\
		&\leq &   \frac{e^{5\eta t }}{(\eta t)^{\frac{1}{3}}} + C \left(\frac{1}{\eta} +2 \right)^{2} \frac{e^{5 \eta t}}{(\eta t)^{\frac{1}{3}}} \leq C\left(1 + \left(\frac{1}{\eta} +2 \right)^{2} \right) \frac{e^{5 \eta t}}{(\eta t)^{\frac{1}{3}}} \\
		&\leq & \frac{C}{\eta^{\frac{1}{3}}} \left( 1 + \left(\frac{1}{\eta} +2 \right)^{2} \right) \frac{e^{5 \eta t}}{t^{\frac{1}{3}}}.  
	\end{eqnarray*} Finally, from now on we set the constant  
	\begin{equation}\label{Const-eta}
	c_\eta =\frac{C}{\eta^{\frac{1}{3}}} \left( 1 + \left(\frac{1}{\eta} +2 \right)^{2} \right)>0,
	\end{equation} and we get the desired estimate. 
	%finpagina8
	\item[2)] We will suppose that there exists $\varepsilon>0$ and $M>0$ such that for  all $\vert x \vert>M$, we have $\ds{\vert K_{\eta}(t,x)\vert \lesssim \frac{1}{\vert x \vert^{2+\varepsilon}}}$ and then we will arrive to a contradiction. Indeed, if we suppose this estimate then  we can prove  that the function  $x K_{\eta}(t,x)$ belongs to the space $L^{1}(\R)$: we write 
	$$ \int_{\R} \vert x K_{\eta}(t,x) \vert dx = \int_{\vert x \vert \leq M} \vert x K_{\eta}(t,x) \vert dx + \int_{\vert x \vert >M} \vert x K_{\eta}(t,x) \vert dx = I_1+I_2.$$
	In order to estimate the term $I_1$, recall that from  point $1)$ of Proposition \ref{Prop1}, we have: for all $t>0$,   $K_{\eta}(t,\cdot)\in L^{1}(\R)$. Thus, we have 
	$$ I_1 \leq M \int_{\vert x \vert \leq M} \vert K_\eta(t,x) \vert dx \leq M \Vert K_\eta(t,\cdot)\Vert_{L^1}<+\infty.$$
	Now, we estimate the term $I_2$ and since we have $\ds{\vert K_{\eta}(t,x)\vert \lesssim \frac{1}{\vert x \vert^{2+\varepsilon}}}$, for all $\vert x \vert>M$, then we can write 
	$$ I_2 \lesssim \int_{\vert x \vert >M} \vert x \vert \frac{1}{\vert x \vert^{2+\varepsilon}} dx \lesssim  \int_{\vert x \vert >M} \frac{dx}{\vert x \vert^{1+\varepsilon}} dx <+\infty.$$
	Thus,  the function $x K_{\eta}(t,x)$ belongs to the space $L^{1}(\R)$ and then by the properties of the Fourier transform we get that  $\partial_\xi \widehat{K_{\eta}}(t,\xi)$  is a continuous function. Moreover,  recall that we have $K_\eta(t,\cdot) \in L^{1}(\R)$ and then $\widehat{K_{\eta}}(t,\xi)$ is also a continuous function and thus,  for all time  $t>0$, we have $\widehat{K_{\eta}}(t,\cdot) \in \mathcal{C}^{1}(\R)$, but  this fact is not possible. Indeed, by identity (\ref{Kernel}), we have $\ds{\widehat{K_{\eta}}(t,\xi)= e^{i\xi^3t} e^{-\eta t \vert \xi \vert^3} e^{\eta t \vert \xi \vert}}$, but observe that the term $\ds{e^{\eta t \vert \xi \vert}}$ is not differentiable at the origin and then $\widehat{K_{\eta}}(t,\cdot)$ cannot belong to the space $\mathcal{C}^{1}(\R)$.  \finpv
\end{enumerate} 
%%%%%%%%%%%%%%%%%%%%%%%%%%%%%%%%%%%%%%%%%%%%%%%%%%%%%%%%%%%%%%%%%%%%%%%%%%%%%%%%%%%%%%%%%%
\section{Pointwise decaying and asymptotic behavior in spacial variable}\label{sec:pointwise-decay-asymptotics}
\subsection{Proof of Theorem \ref{Th-decay}}\label{sec:Th1} 
Let  $\frac{3}{2}<s\leq 2$ fix and let  $u_0 \in  H^{s}(\R)$  be the initial datum and suppose that this function verifies 
\begin{equation}\label{cond-initial-data}
\vert u_0(x)\vert \leq  \frac{c}{1 +\vert x \vert^2}.  
\end{equation}
We start by studying the existence of a local in time solution $u$ of integral  equation (\ref{integral}). 
\subsubsection{Local in time existence}
Let $T>0$  and consider the functional space  
$\ds{Y_T= \left\{u \in \mathcal{S}^{'}([0,T] \times \R) : \sup_{0 < t\leq T} t^{\frac{1}{3}} \Vert (1+\vert \cdot\vert^2) u(t,\cdot)\Vert_{L^{\infty}}<+\infty \right\}}$ and then define the Banach space 
\begin{equation}\label{F-T}
F_T=  Y_T \cap  \mathcal{C}\left([0,T], H^{s}(\R)\right),
\end{equation} doted with the norm
\begin{equation}\label{norme-F-T}
\Vert \cdot \Vert_{F_T}= \sup_{t \in ]0,T] }t^{\frac{1}{3}} \Vert (1+\vert \cdot \vert^2)(\cdot)\Vert_{L^{\infty}(\R)} + \sup_{t \in [0,T]} \Vert \cdot\Vert_{H^{s}(\R)}. 
\end{equation}
Remark that this norm is composed of two terms:  the first term in the right side in (\ref{norme-F-T}) will allow us to study the decay in spatial variable of the solution $u$. In this term we can observe a weight in time variable $t^{\frac{1}{3}}$ which the reason to add this weight  is purely technical and it allows us to carry out the estimates  which we shall need later.  \\
\\ On the other hand, the second term on the right side in (\ref{norme-F-T})  will allow us to study the regularity of solution $u$ and this will be done later in Section \ref{Sec:regularity}.    
\begin{Theoreme}\label{Th-aux-1} There exists a time $T_0 >0$  and a function $u \in F_{T_0}$ which is the unique solution of the integral equation (\ref{integral}). 
\end{Theoreme}	
%%%finpagina9
\pv  We write
\begin{eqnarray}\label{eq11} \nonumber
\Vert u \Vert_{F_T}&=&\left\Vert K_\eta(t,\cdot)\ast u_0 -\frac{1}{2} \int_{0}^{t}K_{\eta}(t-\tau,\cdot)\ast \partial_{x}(u^{2})(\tau,\cdot) d \tau \right\Vert_{F_T} \\
& \leq & \Vert  K_\eta(t,\cdot)\ast u_0 \Vert_{F_T}+ \left\Vert  \frac{1}{2} \int_{0}^{t}K_{\eta}(t-\tau,\cdot)\ast \partial_{x}(u^{2})(\tau,\cdot) d \tau \right\Vert_{F_T},
\end{eqnarray} and  we will estimate each term in the right side. 
\begin{Proposition}\label{Prop3-linear} There exist a constant $C_{1,\eta}>0$ given in the formula (\ref{Const-1-eta}), which only depends on $\eta>0$, such that we have: 
	\begin{equation}\label{eq12}
	\Vert  K_\eta(t,\cdot)\ast u_0 \Vert_{F_T} \leq C_{1,\eta}\, e^{5\eta T}  \left(\Vert (1+\vert \cdot \vert^2)u_0 \Vert_{L^{\infty}} +\Vert u_0 \Vert_{H^s}\right). 
	\end{equation}
\end{Proposition}
\pv By the definition of the quantity $\Vert \cdot \Vert_{F_T} $ given in the  equation (\ref{norme-F-T}) we write 
\begin{equation}\label{eq13}
\Vert  K_\eta(t,\cdot)\ast u_0 \Vert_{F_T} = \sup_{t\in ]0,T]} t^{\frac{1}{3}} \Vert (1+\vert \cdot \vert^2) K_\eta(t,\cdot)\ast u_0 \Vert_{L^{\infty}} + \sup_{t\in [0,T]}  \Vert K_\eta(t,\cdot)\ast u_0 \Vert_{H^{s}}, 
\end{equation}
and we start by  estimating the first term on the right side.  For all $x\in \R$ we write 
\begin{eqnarray}\label{eq14} \nonumber
\vert K_\eta(t,\cdot)\ast u_0 (x)\vert  &\leq& \int_{\R} \vert K_\eta(t,x-y) \vert  \vert u_0 (y)\vert dy \leq \int_{\R} \vert K_\eta(t,x-y) \vert \frac{1+\vert y \vert^2}{1+\vert y \vert^2} \vert u_0(y) \vert dy \\
&\leq & \Vert (1+\vert \cdot \vert^2) u_0 \Vert_{L^{\infty}} \int_{\R} \frac{\vert K_\eta(t,x-y)}{1+\vert y \vert^2} dy.
\end{eqnarray} We need to study the term $\ds{\int_{\R} \frac{\vert K_\eta(t,x-y)}{1+\vert y \vert^2} dy}$. Remark that from point $1)$ of Proposition \ref{Prop1} we have the estimate $ \ds{\vert K_\eta(t,x-y)\vert \leq \frac{c_\eta e^{5\eta t}}{t^{\frac{1}{3}}} \frac{1}{1+\vert x-y\vert^2}}$, and then we can write 
\begin{equation}\label{eq15}
\int_{\R} \frac{\vert K_\eta(t,x-y)\vert}{1+\vert y \vert^2} dy \leq \frac{c_\eta e^{5\eta t}}{t^{\frac{1}{3}}} \int_{\R} \frac{dy}{(1+\vert x-y\vert^2)(1+\vert y \vert^2) },
\end{equation}
where the last term on the right side verifies 
\begin{equation}\label{eq16} 
\int_{\R} \frac{dy}{(1+\vert x-y\vert^2)(1+\vert y \vert^2)} \leq c \frac{1}{1+\vert x \vert^2}.
\end{equation}
%Indeed, for $x\in \mathbb{R}$ fix we write 
%\begin{equation}\label{eq17}
% \int_{\R} \frac{dy}{(1+\vert x-y\vert^2)(1+\vert y \vert^2)} = \int_{\vert y \vert \leq \frac{\vert x \vert }{2}} \frac{dy}{(1+\vert x-y\vert^2)(1+\vert y \vert^2)}+ \int_{\vert y \vert > \frac{\vert x \vert }{2}} \frac{dy}{(1+\vert x-y\vert^2)(1+\vert y \vert^2)},
%\end{equation}  then, for the first term in the right side, since $\vert y \vert \leq \frac{\vert x \vert}{2}$ then we have $\vert x-y \vert \geq \vert x \vert -\vert y \vert \geq  \frac{\vert x \vert }{2}$ and thus we can write 
%$$ \int_{\vert y \vert \leq \frac{\vert x \vert }{2}} \frac{dy}{(1+\vert x-y\vert^2)(1+\vert y \vert^2)} \leq \frac{1}{1+\vert  x \vert^2} \int_{\vert y \vert \leq \frac{\vert x \vert}{2}} \frac{dy}{1+\vert y\vert^2} \leq \frac{1}{1+\vert x \vert^2}\int_{\R} \frac{dy}{1+\vert y\vert^2}\leq \frac{c}{1+\vert x \vert}.$$
%Now, for the second term in the right side in (\ref{eq17}), since $\vert y \vert > \frac{\vert x \vert}{2}$ then we have 
%$$ \int_{\vert y \vert > \frac{\vert x \vert }{2}} \frac{dy}{(1+\vert x-y\vert^2)(1+\vert y \vert^2)} \leq \frac{1}{1+\vert x \vert^2}\int_{\vert y \vert > \frac{\vert x \vert}{2}}\frac{dy}{1+\vert x-y \vert^2}\leq \frac{1}{1+\vert x \vert^2}\int_{\R} \frac{dy}{1+\vert x-y \vert^2}\leq \frac{c}{1+\vert x \vert^2}.$$
%With these estimates we get the estimate given in (\ref{eq16}) and then, getting back to equation (\ref{eq15}) we can write 
%$$ \int_{\R} \frac{\vert K_\eta(t,x-y)}{1+\vert y \vert^2} dy \leq \frac{c_\eta e^{5\eta t}}{t^{\frac{1}{3}}} \frac{1}{1+\vert x \vert^2}.$$
Now, we get back to  (\ref{eq14}) and  we have $\ds{\vert K_\eta(t,\cdot)\ast u_0 (x) \vert \leq \Vert (1+\vert \cdot \vert^2) u_0 \Vert_{L^{\infty}}  \frac{c_\eta e^{5\eta t}}{t^{\frac{1}{3}}} \frac{1}{1+\vert x \vert^2}}$. \\
\\
Thus,  the first term on the right side in (\ref{eq13})  is estimated as follows: 
\begin{equation}\label{eq18}
\sup_{t\in ]0,T]} t^{\frac{1}{3}} \Vert (1+\vert \cdot \vert^2) K_\eta(t,\cdot)\ast u_0 \Vert_{L^{\infty}} \leq c_\eta e^{5\eta T} \Vert (1+\vert \cdot \vert^2)u_0\Vert_{L^{\infty}}. \end{equation}
Now, we go to study the second term on the right side in (\ref{eq13})  and we will prove the following estimate
\begin{equation}\label{eq19}
\sup_{t \in [0,T]}  \Vert K_\eta(t,\cdot)\ast u_0 \Vert_{H^{s}} \leq c  e^{5\eta T} \Vert u_0\Vert_{H^s},
\end{equation}
%finpagina10
where $c>0$ is a numerical constant which does not depend on $\eta>0$.  This estimate relies on the following technical estimate given in Lemma $2.2$, (page 10) of \cite{BorysAlvarez-tesis}: let  $s_1 \in \mathbb{R}$, $\phi\in H^{s_1}(\R)$  and let $s_2\geq 0$. Then, for all $t>0$, we have
\begin{equation}\label{estim-Borys}
\Vert K_{\eta}(t,\cdot)\ast \phi \Vert_{H^{s_1+s_2}} \leq c \frac{e^{5\eta t}}{(\eta t)^{\frac{s_2}{2}}} \Vert \phi \Vert_{H^{s_1}}.   
\end{equation} 
In this estimate  we set $\phi=u_0 \in H^{s}(\R)$, $s_1=s$ and $s_2=0$; and then, for all $0\leq t\leq T$, we get $$  \Vert K_{\eta}(t,\cdot)\ast u_0 \Vert_{H^s} \leq c  e^{5\eta t }\Vert u_0 \Vert_{H^s} \leq c  e^{5\eta T} \Vert u_0 \Vert_{H^s},$$ hence, we have the estimate (\ref{eq19}). Now, by estimates (\ref{eq18}) and (\ref{eq19}) we set the constant $C_{1,\eta}>0$ as 
\begin{equation}\label{Const-1-eta}
C_{1,\eta}= c_\eta+c, 
\end{equation} 
where $c_\eta>0$ is the constant given in the formula (\ref{Const-eta}), and then  we have the estimate given in (\ref{eq12}).  Proposition \ref{Prop3-linear} is proven. \finpv 
Now, we estimate  the second term on the right side in the equation (\ref{eq11}).
\begin{Proposition}\label{Prop3-bilinear} There exists a constant $C_{2,\eta}>0$ given in the  formula (\ref{Const-2-eta}), which  depends only on $\eta>0$, such for all $u\in F_T$ we have 
	\begin{equation}\label{eq20}  
	\left\Vert \frac{1}{2}\int_{0}^{t} K_\eta(t-\tau,\cdot)\ast \partial_{x}(u^2)(\tau,\cdot)d \tau \right\Vert_{F_T} \leq C_{2,\eta}\,  e^{5\eta T} \max(T^{\frac{2}{3}}, T^{\frac{1}{2}}) \Vert u \Vert_{F_T} \Vert u \Vert_{F_T}.
	\end{equation}
\end{Proposition}	
\pv By definition of the norm $\Vert \cdot \Vert_{F_T}$ given in (\ref{norme-F-T}), we write 
\begin{eqnarray}\label{eq21} \nonumber
\left\Vert \frac{1}{2}\int_{0}^{t} K_\eta(t-\tau,\cdot)\ast \partial_{x}(u^2)(\tau,\cdot)d \tau \right\Vert_{F_T} &=& \sup_{t\in [0,T]} t^{\frac{1}{3}} \left\Vert (1+\vert  \cdot\vert^2) \left( \frac{1}{2}\int_{0}^{t} K_\eta(t-\tau,\cdot)\ast \partial_{x}(u^2)(\tau,\cdot)d \tau\right) \right\Vert_{L^{\infty}} \\
& & + \sup_{t\in [0,T]}  \left\Vert  \frac{1}{2}\int_{0}^{t} K_\eta(t-\tau,\cdot)\ast \partial_{x}(u^2)(\tau,\cdot)d \tau \right\Vert_{H^{s}},
\end{eqnarray} and we will estimate each term in the right side. \\
\\
For the first term in (\ref{eq21}), for  all $t\in [0,T]$ we have 
\begin{equation*}
t^{\frac{1}{3}} \left\Vert (1+\vert  \cdot\vert^2) \left( \frac{1}{2}\int_{0}^{t} K_\eta(t-\tau,\cdot)\ast \partial_{x}(u^2)(\tau,\cdot)d \tau\right) \right\Vert_{L^{\infty}} \leq t^{\frac{1}{3}} \int_{0}^{t} \left\Vert (1+\vert \cdot \vert^2) \frac{1}{2}K_\eta(t-\tau, \cdot)\ast \partial_{x}(u^2)(\tau,\cdot)\right\Vert_{L^{\infty}}d\tau,
\end{equation*} and now  we need to  prove the following estimate:

\begin{equation}\label{eq26}
\left\Vert (1+\vert \cdot \vert^2) \frac{1}{2} K_{\eta}(t-\tau,\cdot)\ast \partial_{x}(u^2)(\tau,\cdot) \right\Vert_{L^{\infty}}\leq  c_\eta \frac{e^{5\eta (t-\tau)}}{(t-\tau)^{\frac{1}{3}} \tau^{\frac{1}{3}}}  \Vert u \Vert_{F_T} \Vert u \Vert_{F_T}. 
\end{equation}
Indeed,  we will study first the quantity $ \frac{1}{2}  K_\eta(t-\tau,\cdot)\ast \partial_x (u^2)(\tau,\cdot)(x)$. Remark that 
we have $\frac{1}{2} \partial_x (u^2)=u \partial_x u$ and then for all $x\in \R$ we write 
%%finalpagina11
\begin{eqnarray}\label{eq22} \nonumber
\left\vert \frac{1}{2} K_{\eta}(t-\tau,\cdot)\ast \partial_x (u^2)(\tau,\cdot)(x)\right\vert &\leq & \vert K_{\eta}(t-\tau,\cdot)\ast \left( u(\tau, \cdot) \partial_x u(\tau, \cdot) \right) (x)\vert \\
&\leq & \int_{\R} \vert K_{\eta}(t-\tau, x-y) \vert \vert u(\tau, y)\vert \vert \partial_y u(\tau,y)\vert dy.
\end{eqnarray}
Now, recall that by point $1)$ of  Proposition \ref{Prop1}, we have $\ds{\vert K_\eta(t-\tau,x-y)\vert \leq c_\eta \frac{e^{5\eta (t-\tau)}}{(t-\tau)^{\frac{1}{3}}}\frac{1}{1+\vert x-y \vert^2}}$, and then in the last term above, we can write 
\begin{eqnarray}\label{eq23} 
& &\int_{\R} \vert K_{\eta}(t-\tau, x-y) \vert \vert u(\tau, y)\vert \vert \partial_y u(\tau,y)\vert dy \leq   c_\eta \frac{e^{5\eta (t-\tau)}}{(t-\tau)^{\frac{1}{3}}} \int_{\R} \frac{\vert u(\tau, y) \vert \vert \partial_y u (\tau,u)\vert  }{1+\vert x-y\vert^2} dy \\ \nonumber
&\leq & c_\eta \frac{e^{5\eta (t-\tau)}}{(t-\tau)^{\frac{1}{3}}} \Vert (1+\vert \cdot\vert^2)u(\tau, \cdot)\Vert_{L^{\infty}} \int_{\R} \frac{\vert \partial_y u(\tau, y\vert)}{(1+\vert y \vert^2)(1+\vert x-y\vert^2)}dy \\ 
&\leq & c_\eta \frac{e^{5\eta (t-\tau)}}{(t-\tau)^{\frac{1}{3}}} \underbrace{\Vert (1+\vert \cdot\vert^2)u(\tau, \cdot)\Vert_{L^{\infty}} \Vert \partial_x u(\tau, \cdot)\Vert_{L^{\infty}}}_{(a)}  \underbrace{\int_{\R} \frac{dy}{(1+\vert y \vert^2)(1+\vert x-y\vert^2}}_{(b)},  \label{eq49}
\end{eqnarray}
where we have to study the terms $(a)$ and $(b)$.  For term $(a)$ we have 
\begin{equation}\label{eq24}
(a)\leq  \frac{c}{\tau^{\frac{1}{3}}} \Vert u \Vert_{F_T} \Vert u \Vert_{F_T}.
\end{equation}
Indeed,  recall first that  we have the inclusion $H^{s-1}(\R) \subset L^{\infty}(\R)$ (since $s-1>\frac{1}{2}$). Hence, we can write 
\begin{equation}\label{eq87}
\Vert \partial_{y}u(\tau, \cdot) \Vert_{L^{\infty}} \leq c\Vert \partial_{x}u(\tau, \cdot)\Vert_{H^{s-1}} \leq c\Vert u(\tau,\cdot)\Vert_{H^{s}}. 
\end{equation}  Thus, we have 
\begin{equation*}
(a)\leq \Vert (1+\vert \cdot\vert^2)u(\tau, \cdot)\Vert_{L^{\infty}} \Vert  u(\tau, \cdot)\Vert_{H^{s}}\leq  \frac{c}{\tau^{\frac{1}{3}}} \left( \tau^{\frac{1}{3}} \Vert (1+\vert \cdot\vert^2)u(\tau, \cdot)\Vert_{L^{\infty}}\right) \left(\Vert  u(\tau, \cdot)\Vert_{H^{s}}\right), 
\end{equation*}  
and by definition of the norm $\Vert \cdot \Vert_{F_T}$ given in (\ref{norme-F-T}) we can write the estimate given in (\ref{eq24}).\\
\\
For term $(b)$ in (\ref{eq49}), recall that this  was already estimated at (\ref{eq16}).\\
\\
Then, in the estimate (\ref{eq49}),  by estimates (\ref{eq24}) and (\ref{eq16}) we have
\begin{equation*}
\int_{\R} \vert K_{\eta}(t-\tau, x-y) \vert \vert u(\tau, y)\vert \vert \partial_y u(\tau,y)\vert dy \leq  c_\eta \frac{e^{5\eta (t-\tau)}}{(t-\tau)^{\frac{1}{3}} \tau^{\frac{1}{3}}} \frac{1}{1+\vert x \vert^2} \Vert u \Vert_{F_T} \Vert u \Vert_{F_T}, 
\end{equation*}
and now, we get back to estimate (\ref{eq22}) and we write 
\begin{equation*}
\left\vert \frac{1}{2} K_{\eta}(t-\tau,\cdot)\ast \partial_x (u^2)(\tau,\cdot)(x)\right\vert\leq c_\eta \frac{e^{5\eta (t-\tau)}}{(t-\tau)^{\frac{1}{3}} \tau^{\frac{1}{3}}} \frac{1}{1+\vert x \vert^2} \Vert u \Vert_{F_T} \Vert u \Vert_{F_T}.  
\end{equation*} Thus, we get the estimate (\ref{eq26}).\\
\\
Once we dispose of this estimate,  for all $t\in [0,T]$,  we can write  
\begin{eqnarray*}
	& &  t^{\frac{1}{3}} \int_{0}^{t} \left\Vert (1+\vert \cdot \vert^2) \frac{1}{2}K_\eta(t-\tau, \cdot)\ast \partial_{x}(u^2)(\tau,\cdot)\right\Vert_{L^{\infty}}d\tau  \leq   c_\eta\, t^{\frac{1}{3}}\left(\int_{0}^{t}  e^{5\eta (t-\tau)} \frac{d \tau}{(t-\tau)^{\frac{1}{3}} \tau ^{\frac{1}{3}}} \right) \Vert u \Vert_{F_T} \Vert u \Vert_{F_T} \\
	&\leq & c_\eta \,  t^{\frac{1}{3}}  e^{5 \eta T} \left(\int_{0}^{t} \frac{d \tau}{(t-\tau)^{\frac{1}{3}} \tau^{\frac{1}{3}}} \right) \Vert u \Vert_{F_T} \Vert u \Vert_{F_T} \leq c_\eta \,   T^{\frac{1}{3}} e^{5 \eta } \, \left(  T^{\frac{1}{3}} \right) \Vert u \Vert_{F_T} \Vert u \Vert_{F_T}\\
	&\leq & c_\eta \, e^{5\eta T} \, T^{\frac{2}{3}}  \Vert u \Vert_{F_T} \Vert u \Vert_{F_T}, 
\end{eqnarray*} and then we have 
\begin{equation}\label{eq27}
\sup_{t\in ]0,T]} t^{\frac{1}{3}} \left\Vert (1+\vert  \cdot\vert^2) \left( \frac{1}{2}\int_{0}^{t} K_\eta(t-\tau,\cdot)\ast \partial_{x}(u^2)(\tau,\cdot)d \tau\right) \right\Vert_{L^{\infty}}\leq c_\eta\, e^{5\eta T} \, T^{\frac{2}{3}} \Vert u \Vert_{F_T} \Vert u \Vert_{F_T}. \\
\end{equation}
%FINALPAGINA12

Now, we estimate the second term in identity (\ref{eq21}). For all $t\in [0,T]$, we write 
\begin{eqnarray*}
	& & \left\Vert \int_{0}^{t} K_\eta(t-\tau, \cdot)\ast \partial_x (u^2)(\tau,\cdot) d \tau \right\Vert_{H^{s}} \leq  \int_{0}^{t} \Vert K_{\eta}(t-\tau, \cdot)\ast \partial_x (u^2)(\tau,\cdot)\Vert_{H^{s}} d \tau\\
	&\leq & \int_{0}^{t}\Vert \partial_x (K_\eta(t-\tau, \cdot) \ast u^{2}(\tau, \cdot))\Vert_{H^{s}} d\tau \leq \int_{0}^{t}\Vert  K_\eta(t-\tau, \cdot) \ast u^{2}(\tau, \cdot)\Vert_{H^{s+1}} d\tau.
\end{eqnarray*}
Then, in the estimate (\ref{estim-Borys}) we set now $\phi= (u^2)(\tau,\cdot)$, $s_1=s$ and $s_2=1$; and then we have 
\begin{equation*}
\int_{0}^{t}\Vert  K_\eta(t-\tau, \cdot) \ast u^{2}(\tau, \cdot)\Vert_{H^{s+1}} d\tau \leq \int_{0}^{t} c  \frac{e^{5 \eta (t-\tau)}}{ (\eta(t-\tau))^{\frac{1}{2}}} \Vert u^2(\tau, \cdot)\Vert_{H^{s}} d \tau,
\end{equation*} 
where,  by the product laws in Sobolev spaces and moreover, by  definition of the norm $\Vert \cdot \Vert_{F_T}$ given in (\ref{norme-F-T}), we have 
\begin{eqnarray*}
	& & \int_{0}^{t} c \frac{e^{5 \eta (t-\tau)}}{(\eta (t-\tau))^{\frac{1}{2}}} \Vert u^2(\tau, \cdot)\Vert_{H^{s}} d \tau \leq \int_{0}^{t} c\frac{e^{5 \eta (\eta(t-\tau))}}{(t-\tau)^{\frac{1}{2}}} \Vert u (\tau, \cdot)\Vert^{2}_{H^{s}} d \tau\\
	&\leq & c \frac{e^{5\eta T}}{\eta^{\frac{1}{2}}}  \left(\sup_{\tau \in [0,T]} \Vert u(\tau, \cdot)\Vert_{H^s} \right)\left(\sup_{\tau \in [0,T]} \Vert u(\tau, \cdot)\Vert_{H^s} \right) \int_{0}^{t}\frac{d \tau}{(t-\tau)^{\frac{1}{2}}} \leq c \frac{e^{5\eta T} }{\eta^{\frac{1	}{2}}} T^{\frac{1}{2}} \Vert u \Vert_{F_T} \Vert u \Vert_{F_T}.
\end{eqnarray*}
Thus, we get the estimate
\begin{equation}\label{eq28}
\sup_{t\in [0,T]}   \left\Vert  \frac{1}{2}\int_{0}^{t} K_\eta(t-\tau,\cdot)\ast \partial_{x}(u^2)(\tau,\cdot)d \tau \right\Vert_{H^{s}} \leq c \frac{e^{5\eta T} }{\eta^{\frac{1	}{2}}} T^{\frac{1}{2}} \Vert u \Vert_{F_T} \Vert u \Vert_{F_T}.
\end{equation}
Finally, by estimates (\ref{eq27}) and (\ref{eq28}) we set the constant $C_{2,\eta}>0$ as 
\begin{equation}\label{Const-2-eta}
C_{2,\eta}= c_\eta+ \frac{c}{\eta^{\frac{1}{2}}}, 
\end{equation} 
where $c_\eta>0$ is always the constant given in the formula (\ref{Const-eta}), and the estimate (\ref{eq20}) follows.   Proposition \ref{Prop3-bilinear} in now proven. \finpv
\\	
Once we have the estimates given in Proposition \ref{Prop3-linear} and in Proposition \ref{Prop3-bilinear}, we fix the time $T_0>0$ small enough and by the Picard contraction principle we get a solution $u\in F_{T_0} $ of the integral equation (\ref{integral}). \\
\\
Now, we prove the uniqueness of this solution $u\in F_{T_0}$.  Let $u_1,u_2 \in F_{T_0}$ be two solutions of the equation (\ref{integral}) (associated with the same initial datum $u_0$). We define $v=u_1-u_2$ and we will prove that $v=0$. Indeed, recall first that $v(0,\cdot)=0$ and then $v$ verifies the following integral equation 
\begin{equation*}
v(t,\cdot)=-\frac{1}{2} \int_{0}^{t} K_\eta(t-\tau, \cdot)\ast \left(  \partial_{x} (u^{2}_{1}(\tau,\cdot)- u^{2}_{1}(\tau,\cdot)) \right)d \tau.
\end{equation*}  
Since, $v=u_1-u_2$, we write  $\ds{u^{2}_{1}(\tau,\cdot) -u^{2}_{1}(\tau,\cdot)= v(\tau,\cdot) u_1(\tau,\cdot)+u_2(\tau,\cdot) v(\tau,\cdot)}$, and thus we have 
\begin{equation}\label{v}
v(t,\cdot)=-\frac{1}{2} \int_{0}^{t} K_\eta(t-\tau, \cdot)\ast \left(  \partial_{x} (v(\tau,\cdot) u_1(\tau,\cdot)+u_2(\tau,\cdot) v(\tau,\cdot)) \right)d \tau.
\end{equation} 
%%FINALPAGINA13
In this expression we take the norm $\Vert \cdot \Vert_{F_{T_0}}$ given in (\ref{norme-F-T})  and by Proposition \ref{Prop3-bilinear},  we have

\begin{equation}\label{eq29}
\Vert v \Vert_{F_{T_0}} \leq C_{2,\eta} \max(T^{\frac{2}{3}}_{0}, T^{\frac{1}{2}}_{0}) \Vert v \Vert_{F_{T_0}} \left( \Vert u_1 \Vert_{F_{T_0}}+ \Vert u_2 \Vert_{F_{T_0}} \right). 
\end{equation}
From this estimate,  the identity $v=0$ is deduced as follows: let $0\leq T^{*}\leq T_0$ be the maximal time such that $v=0$ at the interval $[0,T^{*}[$. We will prove  that $T^{*}=T_0$ and by contradiction.
\\
\\
Let us suppose $T^{*}<T_0$. 
Let $T_1 \in ]T^{*}, T_0[$ and for the interval in time $]T^{*},T_1[$, consider the space $F_{(T_1-T^{*})}$ defined in (\ref{F-T}) and endowed with the norm $\Vert \cdot \Vert_{F_{(T_1-T^{*})}}$ given in (\ref{norme-F-T}). By estimate (\ref{eq29}),  we can write 
\begin{equation*}
\Vert v \Vert_{F_{(T_1-T^{*})}} \leq C_{2,\eta} \max\left((T_1-T^{*})^{\frac{2}{3}}, (T_1-T^{*})^{\frac{1}{2}}\right)  \Vert v \Vert_{F_{(T_1-T^{*})}}  \left( \Vert u_1 \Vert_{F_{(T_1-T^{*})}}+ \Vert u_2 \Vert_{F_{(T_1-T^{*})}} \right),
\end{equation*} 
and taking $T_1-T^{*}>0$ small enough, then we have $\ds{\Vert v \Vert_{F_{(T_1-T^{*})}}}=0$ and thus we have $v=0$  in the interval in time $]T^{*},T_1[$, which is a contraction with the definition of time $T^{*}$. Then we have $T^{*}=T$. Theorem \ref{Th-aux-1} is now proven. \finpv
\subsubsection{Global in time existence and decay in spacial variable}
In this section, we prove first that the local in time solution $u\in F_{T_0}$ of the integral equation (\ref{integral})  is extended to the whole interval in time $]0,+\infty[$. Then, we prove the decay in spatial variable given in  the formula (\ref{decay}). 
\begin{Theoreme}\label{Th-aux-2} Let $T_0>0$ be the time given in Theorem \ref{Th-aux-1}. Let the Banach space $(F_{T_0}, \Vert \cdot \Vert_{F_{T_0}})$ given by formulas (\ref{F-T}) and (\ref{norme-F-T}) and let $u\in F_{T_0}$ the solution of the integral equation (\ref{integral}) constructed in Theorem \ref{Th-aux-1}. Then, we have:  
	\begin{enumerate}
		\item[1)]  $u \in \mathcal{C}([0,+\infty[, H^s(\R))$.
		\item[2)]  Moreover, for all time $t>0$, there exists a constant $C=C(t,\eta, u_0, \|u(t)\|_{H^{s}})>0$, which  depends on $t>0$, $\eta >0$,  $u_0$, and the quantity $\| u(t) \|_{H^{s}}$, such that, for all $x\in \R$, the solution $u(t,x)$ verifies the estimate (\ref{decay}). 
		%For all time $T>0$  set the quantity $\ds{C_1(T,u)= \sup_{t \in [0,T]}\Vert u(t,\cdot)\Vert_{H^s} <+\infty}$. Then, there exists a constant  $0<C=C\left(T,\eta, u_0,C_1 (T,u)\right)<+\infty$, which  depends of $T>0$, $\eta >0$, the initial data $u_0$ and the quantity $C_1(T,u)$,   such that we have 
		%\begin{equation}\label{eq50}	
		%\sup_{t \in [0,T]} t^{\frac{1}{3}} \Vert (1+\vert \cdot \vert^2) u(t,\cdot)\Vert_{L^{\infty}} \leq C.  \\
		%\end{equation}
	\end{enumerate}		
\end{Theoreme}	
%In estimate (\ref{eq50}) we can observe  that since the constant $C=C\left(T,\eta, u_0,C_1 (T,u)\right)$ is well-defined for all time $T>0$ then the quantity $\ds{\sup_{t \in [0,T_0]} t^{\frac{1}{3}} \Vert (1+\vert \cdot \vert^2) u(t,\cdot)\Vert_{L^{\infty}}}$ does not explode in a finite time and thus the solution $ u \in \mathcal{C}([0,+\infty[, H^s(\R))$ verifies the estimate (\ref{eq50}) for all time $T>0$. \\
%\\
%Remark also that in estimate (\ref{eq50}) we obtain the decay  in spatial variable of the solution $u$ stated in formula (\ref{decay}) in Theorem \ref{Th-decay}. \\
%\\
\pv \begin{enumerate}
	\item[1)]  Since $u_0 \in H^{s}(\R)$, we get by Theorem  $2$ of the article \cite{ZhaoCui} that  there exists a function $v \in \mathcal{C}([0,+\infty[, H^s(\R))$, which is the unique solution of integral equation (\ref{integral}). But, by definition of the Banach space $F_T$, we have the inclusion $F_T\subset \mathcal{C}([0,T], H^s(\R))$ and then the solution $u \in F_T$ belongs to the space $\mathcal{C}([0,T], H^s(\R))$. Thus, by the uniqueness of solution $v$, we have $u=v$ on the interval of time $[0,T]$ and then 
	$$ \sup_{t\in [0,T] }\Vert u(t,\cdot)\Vert_{H^{s}}= \sup_{t\in [0,T]} \Vert v(t,\cdot)\Vert_{H^{s}}.$$
	In this identity, we can see that $v \in \mathcal{C}([0,+\infty[, H^s(\R))$  and thus, the quantity $\ds{ \sup_{t\in [0,T]} \Vert u(t,\cdot)\Vert_{H^{s}}}$ does not explode in a finite time and thus the solution $u$  extends to the whole interval of time $[0,+\infty[$. Therefore, we have  $u \in \mathcal{C}([0,+\infty[, H^s(\R))$.\\
	\item[2)] In order to prove the property decay of solution $u\in \mathcal{C}([0,+\infty[, H^s(\R))$ given in  the estimate (\ref{decay}), we will prove that the quantity $\ds{\sup_{t \in ]0,T]} t^{\frac{1}{3}} \Vert (1+\vert \cdot \vert^2) u(t,\cdot)\Vert_{L^{\infty}}}$ is well-defined for all time $T>0$. \\
	\\
	%%FINALPAGINA14
	Let $T>0$.  For all $t\in ]0,T]$, we write  
	%%%%%%%%%%%%%%%%%%%%%%%%%%%%%%%%%%%%%%%%%%%%%
	\begin{eqnarray}\label{eq45} \nonumber 
	t^{\frac{1}{3}} \Vert (1+\vert \cdot \vert^2) u(t,\cdot)\Vert_{L^{\infty}}&\leq &  t^{\frac{1}{3}}\left\Vert (1+\vert \cdot \vert^2) \left( K_\eta(t,\cdot)\ast u_0 -\frac{1}{2}\int_{0}^{t}K_\eta(t-\tau, \cdot)\partial_x (u^{2})(\tau, \cdot)d\tau\right)\right\Vert_{L^{\infty}} \\ \nonumber 
	&\leq &  t^{\frac{1}{3}}\left\Vert (1+\vert \cdot \vert^2) \left( K_\eta(t,\cdot)\ast u_0 \right)\right\Vert_{L^{\infty}}\\ \nonumber
	& & +  t^{\frac{1}{3}}\left\Vert (1+\vert \cdot \vert^2) \left(\frac{1}{2} \int_{0}^{t}K_\eta(t-\tau, \cdot)\partial_x (u^{2})(\tau, \cdot)d\tau\right)\right\Vert_{L^{\infty}} \\
	&\leq & I_1+I_2,
	\end{eqnarray} We will study the terms $I_1$ and $I_2$ above. For  term $I_1$, by  Proposition \ref{Prop3-linear}, we have 
	\begin{equation*}
	I_1 \leq t^{\frac{1}{3}} \Vert (1+\vert \cdot \vert^2) K_\eta(t,\cdot)\ast u_0 \Vert_{L^{\infty}} \leq C_{1,\eta}\, e^{5\eta T} \Vert (1+\vert \cdot \vert^2)u_0 \Vert_{L^{\infty}}, 
	\end{equation*} where we set the constant as
	\begin{equation}\label{C0}
	\mathfrak{C}_0(T,\eta,u_0)=C_{1,\eta}\, e^{5\eta T} \Vert (1+\vert \cdot \vert^2)u_0 \Vert_{L^{\infty}}>0,
	\end{equation} and then, we write 
	\begin{equation}\label{eq37}
	I_1 \leq \mathfrak{C}_0(T,\eta,u_0).
	\end{equation}
	Now, we compute the  $I_2$  on the right side of the formula (\ref{eq45}).  We write 
	\begin{eqnarray}\label{eq46}\nonumber
	I_2&\leq & t^{\frac{1}{3}}\left\Vert (1+\vert \cdot \vert^2) \left( \int_{0}^{t}K_\eta(t-\tau, \cdot)\partial_x (u^{2})(\tau, \cdot)d\tau\right)\right\Vert_{L^{\infty}} \\
	&\leq &  t^{\frac{1}{3}} \int_{0}^{t} \underbrace{\frac{1}{2}\left\Vert (1+\vert \cdot \vert^2) \frac{1}{2}K_{\eta}(t-\tau)\ast \partial_x (u^{2})(\tau,\cdot)\right\Vert_{L^{\infty}}}_{(a)} d\tau, 
	\end{eqnarray} and  we will estimate the term $(a)$. 
	Indeed, the first thing to do is to study the quantity $$\ds{\left\vert \frac{1}{2} K_{\eta}(t-\tau,\cdot)\ast \partial_x (u^2)(\tau,\cdot)(x)\right\vert },$$  and  
	by estimates (\ref{eq22}) and (\ref{eq23}). We have 
	\begin{equation}\label{eq51}
	\left\vert \frac{1}{2} K_{\eta}(t-\tau,\cdot)\ast \partial_x (u^2)(\tau,\cdot)(x)\right\vert  \leq c_\eta \frac{e^{5\eta (t-\tau)}}{(t-\tau)^{\frac{1}{3}}} \int_{\R} \frac{\vert u(\tau, y) \vert \vert \partial_y u (\tau,u)\vert  }{1+\vert x-y\vert^2} dy,  
	\end{equation} where the constant $c_\eta>0$ is given in (\ref{Const-eta}), and then we write
	\begin{eqnarray}\label{eq48} \nonumber 
	& & c_\eta \frac{e^{5\eta (t-\tau)}}{(t-\tau)^{\frac{1}{3}}} \int_{\R} \frac{\vert u(\tau, y) \vert \vert \partial_y u (\tau,u)\vert  }{1+\vert x-y\vert^2} dy  \leq  c_\eta \frac{e^{5\eta T}}{(t-\tau)^{\frac{1}{3}}} \int_{\R} \frac{\vert u(\tau, y) \vert \vert \partial_y u (\tau,u)\vert  }{1+\vert x-y\vert^2} dy \\ \nonumber
	&\leq & c_\eta \frac{e^{5\eta T}}{(t-\tau)^{\frac{1}{3}}\tau^{\frac{1}{3}}} \int_{\R} \frac{ \tau^{\frac{1}{3}}(1+\vert y \vert^2) \vert u(\tau, y) \vert \vert \partial_y u (\tau,u)\vert  }{(1+\vert y\vert^2)(1+\vert x-y\vert^2)} dy \\ 
	&\leq & c_\eta \frac{e^{5\eta T}}{(t-\tau)^{\frac{1}{3}}\tau^{\frac{1}{3}}} \left(\tau^{\frac{1}{3}} \Vert (1+\vert \cdot \vert^2)u(\tau, \cdot)\Vert_{L^{\infty}} \right)\underbrace{\left(\Vert \partial_x u(\tau,\cdot)\Vert_{L^{\infty}} \right)}_{(a.1)} \underbrace{ \int_{\R} \frac{dy}{(1+\vert y\vert^2)(1+\vert x-y\vert^2)}}_{(a.2)},  
	\end{eqnarray}
	where we still need to estimate the terms $(a.1)$ and $(a.1)$. For the term $(a.1)$, always with $s-1>\frac{1}{2}$ and thus, we can write  $\ds{ (a.1) \leq \partial_{x}u(\tau,\cdot)\Vert_{H^{s-1}} \leq \Vert u(\tau,\cdot)\Vert_{H^{s}}}$. Now,  by point $1)$ of Theorem \ref{Th-aux-2}, we have $u\in \mathcal{C}([0,+\infty[, H^{s}(\R)$ and then , we get 
	$\ds{ (a.1) \leq \sup_{\tau \in [0,T]}\Vert u(\tau, \cdot)\Vert_{H^{s}}}$. Thus, we set the quantity 
	%finpagina15
	\begin{equation}\label{C(T,u)}
	\mathfrak{C}_1(T,u)=\sup_{\tau \in [0,T]}\Vert u(\tau, \cdot)\Vert_{H^{s}}>0, 
	\end{equation} and we can write
	\begin{equation}\label{a.1}
	(a.1) \leq \mathfrak{C}_1(T,u).
	\end{equation}
	On the other hand, recall that term $(a.2)$ was estimated in  the formula (\ref{eq16}) by 
	$\ds{ (a.2) \leq c\frac{1}{1+\vert x \vert^2}}$. \\
	\\
	In this way,  we substitute estimates (\ref{a.1}) and (\ref{eq16}) in terms $(a.1)$ and $(a.2)$ respectively given in the formula (\ref{eq48}), and we get 
	\begin{eqnarray}\label{eq82} \nonumber
	& & c_\eta \frac{e^{5\eta T}}{(t-\tau)^{\frac{1}{3}}\tau^{\frac{1}{3}}} \left(\tau^{\frac{1}{3}} \Vert (1+\vert \cdot \vert^2)u(\tau, \cdot)\Vert_{L^{\infty}} \right) \left(\Vert \partial_x u(\tau,\cdot)\Vert_{L^{\infty}} \right)  \int_{\R} \frac{dy}{(1+\vert y\vert^2)(1+\vert x-y\vert^2)}\\
	& \leq & c_\eta \frac{e^{5\eta T}}{(t-\tau)^{\frac{1}{3}}\tau^{\frac{1}{3}}} \left(\tau^{\frac{1}{3}} \Vert (1+\vert \cdot \vert^2)u(\tau, \cdot)\Vert_{L^{\infty}} \right) \mathfrak{C}_1(T,u) \frac{1}{1+\vert x \vert^2}.
	\end{eqnarray} 
	Then,  by formulas (\ref{eq51}), (\ref{eq48})  and (\ref{eq82}), we get the following estimate 
	$$ \left\vert \frac{1}{2} K_{\eta}(t-\tau,\cdot)\ast \partial_x (u^2)(\tau,\cdot)(x)\right\vert c_\eta \frac{e^{5\eta (t-\tau)}}{(t-\tau)^{\frac{1}{3}}\tau^{\frac{1}{3}}} \Vert (1+\vert \cdot \vert^2)u(\tau, \cdot)\Vert_{L^{\infty}} \mathfrak{C}_1(T,u) \frac{1}{1+\vert x \vert^2},$$ and by this estimate, for term $(a)$ given in right side of estimate (\ref{eq46}) we can write 
	\begin{eqnarray*}
		(a) &=& \Vert (1+\vert \cdot \vert^2) K_{\eta}(t-\tau)\ast \partial_x (u^{2})(\tau,\cdot)\Vert_{L^{\infty}}  \leq c_\eta \frac{e^{5\eta T}}{(t-\tau)^{\frac{1}{3}}\tau^{\frac{1}{3}}} \left(\tau^{\frac{1}{3}} \Vert (1+\vert \cdot \vert^2)u(\tau, \cdot)\Vert_{L^{\infty}} \right)\mathfrak{C}_1(T,u) \\
		&\leq & c_\eta \frac{e^{5\eta T} \mathfrak{C}_1(T,u)}{(t-\tau)^{\frac{1}{3}}\tau^{\frac{1}{3}}} \left(\tau^{\frac{1}{3}} \Vert (1+\vert \cdot \vert^2)u(\tau, \cdot)\Vert_{L^{\infty}} \right).
	\end{eqnarray*} 
	Now, we get back to estimate  (\ref{eq46}) and we have 
	\begin{eqnarray*}
		I_2 &\leq& c_\eta\, t^{\frac{1}{3}}   e^{5\eta T} \mathfrak{C}_1(T,u) \int_{0}^{t}  \frac{1}{(t-\tau)^{\frac{1}{3}}\tau^{\frac{1}{3}} } \left(\tau^{\frac{1}{3}} \Vert (1+\vert \cdot \vert^2)u(\tau, \cdot)\Vert_{L^{\infty}} \right) d\tau\\
		&\leq & c_\eta  \,T^{\frac{1}{3}} (e^{5\eta T} \mathfrak{C}_1(T,u))  \int_{0}^{t}  \frac{1}{(t-\tau)^{\frac{1}{3}}\tau^{\frac{1}{3}}} \left(\tau^{\frac{1}{3}} \Vert (1+\vert \cdot \vert^2)u(\tau, \cdot)\Vert_{L^{\infty}} \right) d\tau.
	\end{eqnarray*} 
	At this point, with the constant $c_\eta>0$ given in (\ref{Const-eta}) and  the constant $\mathfrak{C}_1(T,u)$ given in (\ref{C(T,u)}), we set the constant 
	\begin{equation}\label{C2}
	\mathfrak {C}_2(T,\eta, u)= c_\eta\,T^{\frac{1}{3}} (e^{5\eta T} \mathfrak{C}_1(T,u))>0,
	\end{equation} 
	and, then we write 
	\begin{equation}\label{eq52}
	I_2 \leq \mathfrak{C}_2(T,\eta,u) \int_{0}^{t}  \frac{1}{(t-\tau)^{\frac{1}{3}}\tau^{\frac{1}{3}}} \left(\tau^{\frac{1}{3}} \Vert (1+\vert \cdot \vert^2)u(\tau, \cdot)\Vert_{L^{\infty}} \right) d\tau.
	\end{equation} With estimates (\ref{eq37}) and (\ref{eq52}), we get back to estimate (\ref{eq45}), and then for all  $t\in [0,T]$, we can write 
	\begin{equation}\label{eq53}
	t^{\frac{1}{3}} \Vert (1+\vert \cdot \vert^2) u(t,\cdot)\Vert_{L^{\infty}} \leq \mathfrak{C}_0(\eta,T,u_0) + \mathfrak{C}_2(\eta,T,u) \int_{0}^{t}  \frac{1}{(t-\tau)^{\frac{1}{3}}\tau^{\frac{1}{3}}} \left(\tau^{\frac{1}{3}} \Vert (1+\vert \cdot \vert^2)u(\tau, \cdot)\Vert_{L^{\infty}} \right) d\tau.
	\end{equation} Now, in order to prove that  quantity $\ds{t^{\frac{1}{3}} \Vert (1+\vert \cdot \vert^2) u(t,\cdot)\Vert_{L^{\infty}}}$  does not explode in a finite time, we will use the following Gr\"onwall's type inequality. For a proof of this  result see Lemma $7.1.2$ of the book \cite{DHenry}.
	%finalpagina16
	%revisado parte del kernel 23-033 
	\begin{Lemme}\label{Lemme-growall-tech} Let $\beta >0$ and $\gamma>0$, such that $\beta +\gamma>1$. Let  $g:[0,T]\longrightarrow [0,+\infty[$  a  function such that, $g$ verifies:
		\begin{enumerate}
			\item[1)] $g \in L^{1}_{loc}([0,T])$, 
			\item[2)] $t^{\gamma-1} g \in L^{1}_{loc}([0,T])$, and
			\item[3)] there exits two constants $a\geq 0$ and $b\geq 0$  such that for almost all  $t\in [0,T]$, we have 
			\begin{equation}\label{eq54}
			g(t)\leq a +b \int_{0}^{t}(t-\tau)^{\beta-1}\tau^{\gamma-1}g(\tau) d\tau, 
			\end{equation}
		\end{enumerate}	
		then:  
		\begin{enumerate}
			\item[a)] 	There exists a continuous and increasing function $\varTheta:[0,+\infty[ \longrightarrow [0,+\infty[$ defined by 
			\begin{equation}\label{omeg} 
			\ds{\varTheta(t)= \sum_{k=0}^{+\infty}c_k\, t^{\sigma k}},
			\end{equation}
			where $\sigma=\beta+\gamma-1>0$ and  where, for the  Gamma function $\Gamma (\cdot)$ the coefficients $c_k>0$  are given by the recurrence formula: 
			\begin{equation*}
			c_0=1, \quad \text{and}\quad \frac{c_{k+1}}{c_k}= \frac{\Gamma(k \sigma +1)}{\Gamma(k\sigma +\beta +\gamma)}, \quad \text{for}\quad k\geq 1. 
			\end{equation*}
			\item[$b)$] 	For all time $t\in [0,T]$, we have 
			\begin{equation}\label{gronwall-tech}
			g(t) \leq a \varTheta( b^{\frac{1}{\sigma}} \,t).
			\end{equation} 
		\end{enumerate}		
	\end{Lemme} 
	In this lemma, we set $\beta =\frac{2}{3}$, $\gamma= \frac{2}{3}$ (where we have $\beta+\gamma>1$) and we set the function $g(t)=t^{\frac{1}{3}}\Vert (1+\vert \cdot \vert^{2})u(t,\cdot)\Vert_{L^{\infty}}$, which verifies the points $1)$, $2)$ and $3)$ (with $\gamma-1=\frac{1}{3}$). 
	%above.
	%revisarconel Yarin
	%In this lemma, we set $\beta =\frac{2}{3}$, $\gamma= \frac{2}{3}$ (where we have $\beta+\gamma>1$) and we set the function $g(t)=t^{\frac{1}{3}}\Vert (1+\vert \cdot \vert^{2})u(t,\cdot)\Vert_{L^{\infty}}$, which verifies the points $1)$, $2)$ and $3)$ above. Indeed, since $t^{\frac{1}{3}}\Vert (1+\vert \cdot \vert^{2})u(t,\cdot)\Vert_{L^{\infty}}$ then this functions verifies the points $1)$ and $2)$  (with $\gamma-1=\frac{1}{3}$). 
	\\
	\\
	%revisado parte del kernel 23-033 
	On the other hand, if for the constant $\mathfrak{C}_0(T,\eta,u_0)>0$ given in (\ref{C0}) and for the constant  $ \mathfrak{C}_2(T,\eta,u)>0$ given in (\ref{C2}), we set the parameters $a=\mathfrak{C}_0(T,\eta,u_0)>0$, $b=\mathfrak{C}_2(T,\eta,u)>0$.  Moreover, if we set the parameters $\beta-1=-\frac{1}{3}$ and $\gamma-1=-\frac{1}{3}$ then,  we can see that  the point $3)$ is verified by  estimate (\ref{eq53}). Also, remark that  since $\beta=\frac{2}{3}$ and $\gamma=\frac{2}{3}$, then we have $\sigma=\beta +\gamma-1=\frac{1}{3}$ and thus $\frac{1}{\sigma}=3$.
	\\
	\\
	%%%%VOLVER A REVISAR
	Then, by estimate (\ref{gronwall-tech}) of  Lemma \ref{Lemme-growall-tech},  for all time $t\in [0,T]$, we have: for $\ds{b^{\frac{1}{\sigma}}= (\mathfrak{C}_{2}(T,\eta,u))^3 >0}$,
	\begin{equation}\label{eq103}
	t^{\frac{1}{3}}\Vert (1+\vert \cdot \vert^{2}) u(t,\cdot)\Vert_{L^{\infty}} \leq 
	\mathfrak{C}_0(T,\eta,u_0) \varTheta\left(b^{\frac{1}{\sigma}} \, t \right) \leq \mathfrak{C}_0(T,\eta,u_0) \varTheta\left(b^{\frac{1}{\sigma}} \, T \right), 
	\end{equation}
	%hastaaqui
	Finally,  we set the constant  
	\begin{equation}\label{Constante-C}
	C(\eta,t,u_0,u)= \frac{\mathfrak{C}_0(T,\eta,u_0) \varTheta\left((\mathfrak{C}_{2}(T,\eta,u))^3 \, T \right)}{t^{\frac{1}{3}}}>0,
	\end{equation}  and then, we have the estimate given in the formula (\ref{decay}).  Theorem \ref{Th-aux-2} is now proven. \finpv
\end{enumerate}	

\subsubsection{Regularity}\label{Sec:regularity}
In order to finish this proof of Theorem \ref{Th-decay} we will prove now that the solution $u$ of the  equation   is smooth enough is spatial variable.
\begin{Proposition}\label{Prop-smooth} Let $\frac{3}{2}<s\leq 2$ and let $u\in \mathcal{C}([0,+\infty[, H^s(\R))$  be the solution of the integral equation (\ref{integral}) given by point $1)$  of Theorem  \ref{Th-aux-2}. Then, we have $u \in  \mathcal{C}(]0,+\infty[, \mathcal{C}^{\infty}(\R))$.
\end{Proposition}
%finalpagina17	
\pv Recall that by hypothesis on the initial datum $u_0$ given in (\ref*{cond-initial-data}), we have $u_0 \in H^{s}(\R)$ for $\frac{3}{2}<s\leq 2$ and then by Theorem $1$ of the article \cite{ZhaoCui} the solution $u\in \mathcal{C}([0,+\infty[, H^s(\R))$ verifies 
\begin{equation}\label{eq55}
u\in \mathcal{C}\left([0,+\infty[,  \bigcap_{\alpha \geq 0} H^{\alpha}(\R)\right).
\end{equation} 
With this information, we easily deduce the property $\ds{u \in  \mathcal{C}(]0,+\infty[, \mathcal{C}^{\infty}(\R))}$. Indeed, we will prove that for all $k\in \mathbb{N}$, the function $\ds{\partial^{n}_{x}u(t,\cdot)}$ is a H\"older continuous function  on $\R$. Let  $n\in \mathbb{N}$ fix. Then, for $\frac{1}{2}<s_1<\frac{3}{2}$ we set   $\alpha=n+s_1$  and by (\ref{eq55}), we have $\ds{\partial^{n}_{x}u(t,\cdot) \in H^{s_1}(\R)}$.\\
\\
On the other hand, recall that we have the identification  $H^{s_1}(\R)=B^{s_1}_{2,2}(\R)$ (where $B^{s_1}_{2,2}(\R)$ denotes a Besov space \cite{BaChDan}) and moreover we have the inclusion $B^{s_1}_{2,2}(\R) \subset B^{s_1-\frac{1}{2}}_{\infty,\infty }(\R) \subset \dot{B}^{s_1-\frac{1}{2}}_{\infty,\infty }(\R)$. \\
\\
Then, we have $\partial^{n}_{x}u(t,\cdot) \in \dot{B}^{s_1-\frac{1}{2}}_{\infty,\infty }(\R)$.  But, since $\frac{1}{2}<s_1<\frac{3}{2}$, then we have $ 0<s_1-\frac{1}{2}<1$ and thus  $\partial^{n}_{x}u(t,\cdot)$ is a $\beta$- H\"older continuous function with $\beta =s_1-\frac{1}{2}$. 
Theorem \ref{Th-decay} is now proven.  \finpv 
%%revisado
%findelarevisión
%%%%%%%%%%%%%%%%%%%%%%%%%%%%%%%%%%%%%%%%%%%%%%%%%%%%%%%%%%%%%%%%%%%%%%%%%%%%%%%%%%%%%%%%%%
\subsection{Proof of Theorem \ref{Th:asymptotics}}\label{sec:Th2} 
Let $\frac{3}{2} <s \leq 2$ fix, let $u_0 \in H^s(\R)$ be the initial datum and suppose that this function  verifies the following decay properties: for $\varepsilon>0$,  
\begin{equation}\label{cond-initial-data-2}
\vert u_0(x)\vert \leq  \frac{c}{1 +\vert x \vert^{2+\varepsilon}}\quad \text{and} \quad \left\vert \frac{d}{dx}u_0(x)\right\vert \leq \frac{c}{1+\vert x \vert^2}. 
\end{equation}  Let $u \in \mathcal{C}(]0,+\infty[, \mathcal{C}^{\infty}(\R))$ be the solution of equation (\ref{eq:f})  associated with the initial datum $u_0$ above and given by Theorem \ref{Th-decay}.  In order to prove the asymptotic profile of $u(t,x)$ given in formula (\ref{asymptotic}), we write the solution $u(t,x)$ as the integral formulation given in (\ref{integral}) and will study each term on the right-hand  side of the equation  (\ref{integral}). \\
\\
For the first term:  $\ds{K_\eta(t,\cdot)\ast u_0(x)}$, we will prove  the following asymptotic development when $\vert x \vert\longrightarrow +\infty$: 
\begin{equation}\label{eq56}
K_\eta(t,\cdot)\ast u_0(x)= K_{\eta}(t,x)\left(\int_{\R} u_0(y)dy\right) +\, o(t)\left(\frac{1}{\vert x \vert^{2}}\right).
\end{equation} 
Indeed, for  all $t>0$ and  $x\in \R$ we write: 
\begin{eqnarray}\label{eq68} \nonumber 
K_\eta(t,\cdot)\ast u_0(x)&=&\int_{\R} K_\eta(t,x-y)u_0(y)dy = \int_{\R} K_\eta(t,x-y)u_0(y)dy + K_\eta(t,x) \left(\int_{\R}u_0(y)dy\right) \\ \nonumber
& & - K_\eta(t,x) \left(\int_{\R}u_0(y)dy\right) \\ \nonumber
&=& K_\eta(t,x) \left(\int_{\R}u_0(y)dy\right) + \underbrace{\int_{\R} K_\eta(t,x-y)u_0(y)dy}_{(a)}  - \underbrace{K_\eta(t,x) \left(\int_{\R}u_0(y)dy\right)}_{(b)}.
\end{eqnarray}
Now, in expression $(a)$ and expression $(b)$ above, first we cut each integral in two parts: 
\begin{equation}\label{eq71}
 \int_{\mathbb{R}}(\cdot)dy=   + \int_{\vert y \vert <\frac{\vert x\vert}{2}} (\cdot)dy + \int_{\vert y \vert >\frac{\vert x\vert}{2}} (\cdot)dy,
\end{equation}
and then we arrange the terms in order to write 
\begin{eqnarray}\label{eq57} \nonumber
(a)+(b)&=& \int_{\vert y \vert<\frac{\vert x \vert}{2}} \left( K_{\eta}(t,x-y)-K_\eta(t,x) \right)u_0(y)	dy + \int_{\vert y \vert >\frac{\vert x \vert}{2}} K_\eta(t,x-y)u_0(y)dy\\ \nonumber
& & - K_\eta(t,x) \left(\int_{\vert y \vert >\frac{\vert x \vert}{2}} u_0(y)dy\right)\\
&=& I_1+I_2+I_3,
\end{eqnarray} and now, in order to prove identity (\ref{eq56}) we must  prove  that 
\begin{equation}\label{eq66}
I_1+I_2+I_3 = \, o(t)\left(\frac{1}{\vert x \vert^2}\right),\quad \text{when}\quad  \vert x \vert \longrightarrow +\infty.
\end{equation}
In order to study the term $I_1$ in identity (\ref{eq57}) we  need the following technical result.
\begin{Lemme}\label{estim-K-2} Let $t>0$ and let $K_\eta(t,\cdot)$ be the kernel given in (\ref{Kernel}). Then, we have $K_\eta(t,\cdot)\in \mathcal{C}^{1}(\R)$  moreover, there exists a constant $C_\eta>0$, which only depends on $\eta>0$, such that  we have: 
\begin{enumerate}
	\item[1)]  for all $x \neq 0$, $\ds{\vert \partial_x K_\eta(t,x)\vert \leq C_\eta \frac{e^{6\eta t}}{\vert x \vert ^3}}$.
	\item[2)]  $\ds{\vert \partial_x K_\eta(t,x)\vert  \leq C_\eta \frac{e^{6\eta t }}{t^{\frac{2}{3}}} \frac{1}{1+\vert x \vert^3}}$.
\end{enumerate}		
\end{Lemme} 
The proof of this lemma follows essentially the same lines of the proof of point $1)$ of Proposition \ref{Prop1} and  we will postpone this proof for the appendix.  Thus, since $K_\eta(t,\cdot)\in \mathcal{C}^{1}(\R)$  then by Taylor expansion of the first order, for  $\theta=\alpha(x-y)+(1-\alpha)x=x-\alpha y$ and for some $\alpha \in ]0,1[$,  we can write:  
\begin{equation}\label{eq74}
K_\eta(t,x-y)-K_\eta(t,x)=-y \partial_{x} K_{\eta}(t,\theta),
\end{equation}  and then we have  
\begin{equation}\label{eq59}
I_1 \leq  \int_{\vert y \vert \leq \frac{\vert x \vert}{2}} \left\vert K_\eta(t,x-y)- K_\eta(t,x)\right\vert \vert u_0(y)\vert dy \leq \int_{\vert y \vert \leq \frac{\vert x \vert}{2}} \vert y \partial_x K_\eta(t,\theta) \vert \vert u_0(y)\vert dy. 
\end{equation} We estimate now the last term on the right-hand side. Recall first that by point $1)$ of Lemma \ref{estim-K-2} we can write $\ds{\vert  \partial_x K_\eta(t,\theta)\vert }\leq C_\eta \frac{e^{6\eta t}}{\vert \theta \vert^3}$, but since we have $\theta =x-\alpha y$ (with $\alpha \in ]0,1[$)  then we can write $\vert \theta \vert \geq \vert x\vert - \alpha \vert y \vert \geq \vert x \vert -\vert y \vert$ and moreover,  since  we have $\vert y \vert <\frac{\vert x \vert}{2}$ then we write $\vert x \vert -\vert y \vert \geq \frac{\vert x \vert}{2}$, and thus we get $\vert \theta \vert \geq \frac{\vert x \vert}{2}$. Then we have 
\begin{equation}\label{eq75}
\ds{\vert  \partial_x K_\eta(t,\theta)\vert }\leq C_\eta \frac{e^{6\eta t}}{\vert x \vert^3},
\end{equation} and getting back to estimate (\ref{eq59}) we get 
\begin{equation}\label{eq60}
\int_{\vert y \vert \leq \frac{\vert x \vert}{2}} \vert y \partial_x K_\eta(t,\theta) \vert \vert u_0(y)\vert d y \leq C_\eta \frac{e^{6\eta t}}{\vert x \vert^3} \int_{\vert y \vert <\frac{\vert x \vert}{2}}\vert y \vert \vert u_0 (y)\vert dy   \leq C_\eta \frac{e^{6\eta t}}{\vert x \vert^3}  \int_{\R} \vert y \vert \vert u_0 (y)\vert dy, 
\end{equation} where, since the initial datum $u_0$ verifies $\ds{\vert u_0(y)\vert \leq \frac{c}{1+\vert y \vert^{2+\varepsilon}}}$ (with $\varepsilon>0$) then the last term on right-hand side converges. Thus, by estimates (\ref{eq59}) and (\ref{eq60}) we have 
\begin{equation}\label{eq106}
I_1 \leq  \left( C_\eta \, e^{6\eta t}\Vert \,  \vert \cdot \vert u_0 \Vert_{L^1} \right) \frac{1}{\vert x \vert^3},
\end{equation}
and then 
\begin{equation}\label{eq61}
I_1 = o(t)\left(\frac{1}{\vert x \vert^2}\right), \quad \text{when}\quad \vert x \vert \longrightarrow+\infty. 
\end{equation}
Now, for  term $I_2$ in the identity (\ref{eq57}) we write 
\begin{equation}\label{eq62}
I_2 \leq \int_{\vert y \vert >\frac{\vert x \vert}{2}} \vert K_\eta(y,x-y) \vert \vert u_0(y)\vert dy,
\end{equation} and in order to study this term, we have the following estimates: remark that  by point $1$  of Proposition \ref{Prop1} we have \begin{equation}\label{eq64}
\vert K_\eta(t,x-y) \vert \leq c_\eta \frac{e^{5 \eta t}}{t^{\frac{1}{3}}} \frac{1}{1+\vert x-y \vert^2},
\end{equation} 
whence, we get 
\begin{equation}\label{eq79}
\Vert K_\eta (t,\cdot)\Vert_{L^1} \leq c_\eta \frac{e^{5 \eta t}}{t^{\frac{1}{3}}}.
\end{equation}
On the other hand, always by the fact that the initial datum $u_0$ verifies $\ds{\vert u_0(y)\vert \leq \frac{c}{1+\vert y \vert^{2+\varepsilon}}}$ and moreover, since in the term $I_2$ we have $\vert y \vert > \frac{\vert x \vert}{2}$ then, for $\vert x \vert$ large enough, we get 
\begin{equation}\label{eq83}
 \vert u_0 (y)\vert \leq \frac{c}{1+\vert y \vert^{2+\varepsilon}} \leq \frac{c}{\vert y \vert^{2+\varepsilon}} \leq \frac{c}{\vert x \vert^{2+\varepsilon}}.
\end{equation}
With  estimates (\ref{eq79}) and (\ref{eq83}) at hand,  we get back to the formula (\ref{eq62}) and we write 
\begin{eqnarray}\label{eq107} \nonumber
\int_{\vert y \vert >\frac{\vert x \vert}{2}} \vert K_\eta(y,x-y) \vert \vert u_0(y)\vert dy & \leq & \frac{c}{\vert x \vert^{2+\varepsilon}} \int_{\vert y \vert > \frac{\vert x \vert}{2}} \vert K_\eta(t,x-y) \vert dy \leq \frac{c}{\vert x \vert^{2+\varepsilon}} \Vert K_{\eta}(t,\cdot)\Vert_{L^1} \\
&  \leq & \frac{c}{\vert x \vert^{2+\varepsilon}} \left(c_\eta \frac{e^{5 \eta t}}{t^{\frac{1}{3}}} \right), 
\end{eqnarray} and by this estimate and estimate (\ref{eq62})   we have: 
\begin{equation}\label{eq63}
I_2 = \, o(t)\left(\frac{1}{\vert x \vert^2}\right), \quad \text{when}\quad \vert x \vert\longrightarrow +\infty. 
\end{equation}
We study now  the term $I_3$ in the identity (\ref{eq57}).  By the estimate (\ref{eq64}) and for $\vert x \vert$ large enough we can write 
\begin{equation}
I_3\leq \vert K_\eta(t,x) \vert \left(\int_{\vert y \vert >\frac{\vert x \vert}{2}}  \vert u_0(y)\vert dy\right) \leq  c_\eta \frac{e^{5 \eta t}}{t^{\frac{1}{3}}} \frac{1}{\vert x \vert^2}  \left( \int_{\vert y \vert >\frac{\vert x \vert}{2}} \vert  u_0(y)\vert dy\right), 
\end{equation} but, recall that since we have $\ds{\vert u_0(y)\vert \leq \frac{c}{1+\vert y \vert^{2+\varepsilon}}}$ then we get $u_0 \in L^{1}(\R)$, and thus  
we have $$\ds{\lim_{\vert x \vert \longrightarrow +\infty} \left( \int_{\vert y \vert >\frac{\vert x \vert}{2}} \vert  u_0(y)\vert dy\right)=0}.$$ Then we can write
\begin{equation}\label{eq65}
I_3 = \, o(t) \left(\frac{1}{\vert x \vert^2}\right), \quad \text{when} \quad \vert x \vert \longrightarrow +\infty. 
\end{equation} Finally, by the estimates (\ref{eq61}), (\ref{eq63}) and (\ref{eq65})  we get the estimate (\ref{eq66}). \\
\\
Now, for the second term on the right-hand  side in the integral equation (\ref{integral}): $\ds{\frac{1}{2} \int_{0}^{t}K_{\eta}(t-\tau,\cdot)\ast \partial_{x}(u^{2})(\tau,\cdot)(x)d \tau}$, we will prove the following asymptotic profile when $\vert x \vert  \longrightarrow +\infty$: 
\begin{equation}\label{eq69}
\frac{1}{2} \int_{0}^{t}K_{\eta}(t-\tau,\cdot)\ast \partial_{x}(u^{2})(\tau,\cdot)(x)d \tau=\int_{0}^{t} K_\eta(t-\tau, x) \left( \int_{\R}u(\tau,y) \partial_y u(\tau,y) dy \right) d \tau +\, o(t)\left(\frac{1}{\vert x \vert^2}\right). 
\end{equation} Indeed, for all $x\in \R$ we write 
\begin{eqnarray} \label{eq70}  \nonumber 
 \frac{1}{2} \int_{0}^{t} K_{\eta}(t-\tau,\cdot)\ast \partial_{x}(u^{2})(\tau,\cdot)(x)d \tau &=& \int_{0}^{t} K_{\eta_{\eta}}(t-\tau, \cdot) \ast \left( u \, \partial_ x  u (\tau,\cdot) \right)(x) d\tau \\ 
& =& \int_{0}^{t} \underbrace{\int_{\R} K_\eta(t-\tau, x-y) u(\tau, y )\partial_{y}u(\tau, y) d y}_{(c)}   \, d\tau,  
\end{eqnarray} then, in order to study the term $(c)$, following  the same computations done  in the formulas (\ref{eq68}), (\ref{eq71}) and (\ref{eq57}) we write 
\begin{eqnarray*}
(c)&=& K_{\eta}(t-\tau, x) \left( \int_{\R} u(\tau, y ) \partial_y u (\tau, y) dy \right)d \tau \\
& & +\int_{\vert y \vert < \frac{\vert x \vert}{2}} \left(K_{\eta}(t-\tau, x-y)- K_\eta(t-\tau, x)\right) \left( u(\tau, y ) \partial_y u (\tau, y)\right) dy\, d\tau  \\
& & +  \int_{\vert y \vert >\frac{\vert x \vert}{2}} K_\eta(t-\tau, x-y)\left( u(\tau, y ) \partial_y u (\tau, y) \right) dy \, d\tau-  K_{\eta} (t-\tau, x) \left( \int_{\vert y \vert > \frac{\vert x \vert}{2}} u(\tau, y ) \partial_y u (\tau, y) dy \right) d\tau,
\end{eqnarray*} and getting back to the identity (\ref{eq70}) we have:
\begin{eqnarray}\label{eq72} \nonumber
 & & \frac{1}{2} \int_{0}^{t} K_{\eta}(t-\tau,\cdot)\ast \partial_{x}(u^{2})(\tau,\cdot)(x)d \tau =  \int_{0}^{t} K_{\eta}(t-\tau, x) \left( \int_{\R} u(\tau, y ) \partial_y u (\tau, y) dy \right)d \tau  \\ \nonumber 
& & +  \underbrace{\int_{0}^{t} \int_{\vert y \vert < \frac{\vert x \vert}{2}} \left(K_{\eta}(t-\tau, x-y)- K_\eta(t-\tau, x)\right) \left( u(\tau, y ) \partial_y u (\tau, y)\right) dy\, d\tau}_{I_a}  \\ \nonumber
& & + \underbrace{\int_{0}^{t} \int_{\vert y \vert >\frac{\vert x \vert}{2}} K_\eta(t-\tau, x-y)\left( u(\tau, y ) \partial_y u (\tau, y) \right) dy \, d\tau}_{I_b}  \\
& &   - \underbrace{\int_{0}^{t} K_{\eta} (t-\tau, x) \left( \int_{\vert y \vert > \frac{\vert x \vert}{2}} u(\tau, y ) \partial_y u (\tau, y) dy \right) d\tau}_{I_c}.
 \end{eqnarray} Thus, in order to obtain the asymptotic profile given in (\ref{eq69}), we must prove the following estimate:
\begin{equation}\label{eq73}
I_a + I_b + I_c = o(t)\left( \frac{1}{\vert x \vert^{2}}\right),\quad \text{when}\quad \vert x \vert \longrightarrow +\infty. 
\end{equation}
For the term $I_a$,  by the estimates (\ref{eq74}) and (\ref{eq75}) we can write 
\begin{eqnarray}\label{eq76} \nonumber
I_a &\leq & \int_{0}^{t}\int_{\vert y \vert < \frac{\vert x \vert}{2}}  \left\vert K_{\eta}(t-\tau, x-y)- K_\eta(t-\tau, x)\right\vert \vert y \vert \left\vert  u(\tau, y ) \partial_y u (\tau, y)\right\vert  dy \, d\tau \\
&\leq & \int_{0}^{t} \left( C_\eta  \frac{e^{6 \eta (t-\tau)}}{\vert x \vert^{3}}  \int_{\R}   \vert y \vert \left\vert  u(\tau, y ) \partial_y u (\tau, y)\right\vert  dy \right) d\tau  \leq C_\eta \frac{e^{6 \eta t}}{\vert x \vert^{3}} \int_{0}^{t} \int_{\R} \vert y \vert   \left\vert  u(\tau, y ) \partial_y u (\tau, y)\right\vert  dy \, d\tau, 
\end{eqnarray} where, in order to estimate the last term on the right-hand side we have the following technical result.
\begin{Lemme}\label{Prop-tech-u} Since the initial data  $u_0$ verifies $\ds{\left\vert \frac{d}{dx}u_0 (x) \right\vert  \leq \frac{c}{1+\vert x \vert^2}}$ then  there exists a constant $0<C^{*}=C^{*}(t,\eta,u_0 ,\|u\|_{H^{s}})<+\infty$, which depends on $t>0$,  $\eta>0$, the initial data $u_0$ and the  solution $u$, such that for all time $\tau \in [0,t]$ and for all $y\in \mathbb{R}$ we have
\begin{equation}\label{eq77}
   \left\vert  u(\tau, y ) \partial_y u (\tau, y)\right\vert  \leq \frac{C^{*}}{\tau^{\frac{2}{3}} (1+\vert y \vert^4)}. 
\end{equation}
\end{Lemme} 
\pv  The first thing to do is to prove that the function $\partial_y u(\tau, y)$ verifies the following estimate: 
\begin{equation}\label{eq84}
 \ds{\vert \partial_y u(\tau,y) \vert \leq \frac{C^{*}_{1}}{\tau^{\frac{1}{3}} (1+\vert y \vert^2)}},
 \end{equation} where $C^{*}_{1}>0$ is a constant which does not depend on the variable $y$. For this we write the solution $u$ as the integral equation (\ref{integral}), then, in each side of this identity (\ref{integral})  we derive respect to the spacial variable $y$  and we have 
$$ \partial_y u(\tau,y)= K_\eta(\tau,\cdot)\ast (\partial_y u_0)(y) -\frac{1}{2} \int_{0}^{\tau} (\partial_y K_{\eta}(\tau-\zeta,\cdot)) \ast \partial_{y}(u^{2})(\zeta,\cdot)(y)d \zeta = I_1+I_2,$$ and now we must study the terms $I_1$ and $I_2$ above. \\
\\
In order to study  term $I_1$, recall that by the second estimate in formula (\ref{cond-initial-data-2}) the initial datum $u_0$ verifies $\ds{\vert \partial_y u_0(y)\vert \leq \frac{c}{1+\vert y \vert^2}}$ and then, in the estimate (\ref{eq18}) we can substitute the function $u_0$ by the function $\partial_y u_0$ and thus  we can write 
\begin{equation}\label{eq85}
\vert I_1\vert  \leq  \vert K_{\eta}(\tau,\cdot) \ast  (\partial_y u_0)(y) \vert\leq  c_\eta \frac{e^{5 \eta \tau}}{\tau^{\frac{1}{3}}} \frac{\Vert 1+ \vert \cdot \vert^2 \partial_y u_0\Vert_{L^{\infty}}}{1+\vert y \vert^2} \leq c_\eta \frac{e^{5 \eta t}}{\tau^{\frac{1}{3}}} \frac{\Vert 1+ \vert \cdot \vert^2 \partial_y u_0\Vert_{L^{\infty}}}{1+\vert y \vert^2},
\end{equation} 
We study now the term $I_2$ and for this we write 
\begin{equation}\label{eq86}
\vert I_2 \vert  \leq \left\vert \frac{1}{2} \int_{0}^{\tau} (\partial_y K_{\eta}(\tau-\zeta,\cdot)) \ast \partial_{y}(u^{2})(\zeta,\cdot)(y)d \zeta \right\vert 
 \leq  \int_{0}^{\tau} \int_{\R}\underbrace{\left\vert \partial_y K_\eta(\tau-\zeta, y-z) \right\vert}_{(a)} \underbrace{ \left\vert \partial_z (u^2)(\zeta, z )\right\vert}_{(b)}  dz \, d \zeta,
\end{equation} where we still need to study the terms $(a)$ and $(b)$. For the term $(a)$ recall that by  point $2)$ of Lemma \ref{estim-K-2}  we have 
\begin{equation}\label{eq89}
\left\vert \partial_y K_\eta(\tau-\zeta, y-z) \right\vert \leq C_\eta \frac{e^{6 \eta (\tau-\zeta)}}{(\tau-\zeta)^{\frac{2}{3}}}\frac{1}{1+\vert y-z \vert^{3}}.
\end{equation}
On the other hand, for the term $(b)$ we have the following estimates
\begin{eqnarray}\label{eq88}\nonumber
\vert \partial_z (u^2)(\zeta, z)\vert &= &2 \vert u(\zeta, z )\vert \vert \partial_z u(\zeta,z)\vert = 2 \frac{(1+\vert z\vert^2)\vert u(\zeta,z \vert)\vert \partial_z u(\zeta,z)\vert }{1+\vert z \vert^2} = 2 \frac{ \zeta^{\frac{1}{3}}(1+\vert z\vert^2)\vert u(\zeta,z \vert)\vert \partial_z u(\zeta,z)\vert }{\zeta^{\frac{1}{3}}(1+\vert z \vert^2)}\\
& \leq & \left( \sup_{0<\zeta <t} \zeta^{\frac{1}{3}} \Vert (1+\vert \cdot \vert^2) u(\zeta,\cdot)\Vert_{L^{\infty}}\right) \left( \sup_{0<\zeta<t} \Vert \partial_z u(\zeta, \cdot)\Vert_{L^{\infty}} \right)\frac{1}{\zeta^{\frac{1}{3}}(1+\vert z \vert^2)},
\end{eqnarray} but, using the quantity $\Vert u \Vert_{F_t}$ (where the norm $\Vert \cdot \Vert_{F_t}$ is  given in the formula (\ref{norme-F-T})) we can write $$  \sup_{0<\zeta <t} \zeta^{\frac{1}{3}} \Vert (1+\vert \cdot \vert^2) u(\zeta,\cdot)\Vert_{L^{\infty}} \leq \Vert u \Vert_{F_t},$$ and moreover, by the estimate (\ref{eq87}) we can write $\ds{\sup_{0<\zeta<t} \Vert \partial_z u(\zeta, \cdot)\Vert_{L^{\infty}} \leq \Vert u \Vert_{F_t}}$,  and thus, getting back to the estimate (\ref{eq88}) we get 
\begin{equation}\label{eq90}
\vert \partial_z (u^2)(\zeta, z)\vert \leq \Vert u \Vert^{2}_{F_t} \frac{1}{\zeta^{\frac{1}{3}}(1+\vert z \vert^2)}. 
\end{equation}
Once we dispose of the estimates (\ref{eq89}) and (\ref{eq90}), we get back to estimate (\ref{eq86}) and then we write 
\begin{eqnarray}\label{eq102} \nonumber
\vert I_2 \vert & \leq & \int_{0}^{t}\int_{\R} \left( C_\eta \frac{e^{6 \eta (\tau-\zeta)}}{(\tau-\zeta)^{\frac{2}{3}}}\frac{1}{1+\vert y-z \vert^{3}} \right) \left(\Vert u \Vert^{2}_{F_t} \frac{1}{\zeta^{\frac{1}{3}}(1+\vert z \vert^2)} \right)dz \, d\zeta \\ \nonumber
& \leq & C_\eta e^{6 \eta \tau} \Vert u \Vert^{2}_{F_t} \left(\int_{0}^{t} \frac{d \zeta }{((\tau -\zeta)^{\frac{2}{3}}) \zeta^{\frac{1}{3}}}\right) \left( \int_{\R} \frac{d z }{(1+ \vert y-z\vert^3)(1+\vert z \vert^2)}  \right) \\ \nonumber
&\leq & C_\eta e^{6 \eta \tau}  \left( \int_{\R} \frac{d z }{(1+ \vert y-z\vert^3)(1+\vert z \vert^2)}  \right) \leq C_\eta e^{6 \eta \tau}  \left( \int_{\R} \frac{d z }{(1+\vert y-z\vert^2)(1+ \vert z \vert^2)} \right) \\
& \leq &  C_\eta e^{6 \eta \tau} \frac{1}{1+\vert y \vert^2} \leq  C_\eta \tau^{\frac{1}{3}} e^{6 \eta \tau} \frac{1}{\tau^{\frac{1}{3}}(1+\vert y \vert^2)} \leq  C_\eta t^{\frac{1}{3}} e^{6 \eta t} \frac{1}{\tau^{\frac{1}{3}}(1+\vert y \vert^2)}. 
\end{eqnarray}
By the estimates (\ref{eq85}) and (\ref{eq102}), we set the constant $C^{*}_{1}$ as $\ds{C^{*}_{1} = \max\left( c_\eta e^{5 \eta t} \Vert (1+\vert \cdot\vert^{2})\partial_y u_0\Vert_{L^{\infty}} , C_\eta t^{\frac{1}{3}} e^{6 \eta t}\right)}>0$, and then we can write the estimate  (\ref{eq77}). \\
\\
Finally, recall that the by estimate (\ref{eq103}) we can write $\ds{\vert u(\tau, y)\vert \leq \frac{\mathfrak{C}_0(t,\eta,u_0) \varTheta\left(b^{\frac{1}{\sigma}} \, t \right)}{\tau^{\frac{1}{3}}(1+\vert y \vert^2)}}$. Thus, we set the constant $C^{*}$ as $C^{*}= \max\left(\mathfrak{C}_0(t,\eta,u_0) \varTheta\left(b^{\frac{1}{\sigma}} \, t \right), C^{*}_{1} \right)>0$,  and then by the estimate above and the estimate (\ref{eq77}) we get the desired estimate (\ref{eq77}).  \finpv 
Thus, getting back to the estimate (\ref{eq76}),  for $\vert x \vert$ large enough  we can write 
$$ I_a \leq  C_\eta \frac{e^{6 \eta t}}{\vert x \vert^{3}} \left( \int_{0}^{t}\int_{\R} \frac{C^{*}}{\tau^{\frac{2}{3}}(1+\vert y \vert^4)} dy\, d\tau \right) \leq  C_\eta \frac{e^{6 \eta t}}{\vert x \vert^{3}} \left( C^{*}\, \left( \int_{0}^{t} \frac{d \tau}{\tau^{\frac{2}{3}}}  \right)\left( \int_{\R} \frac{\vert y \vert}{1+\vert y \vert^4} dy \right)\right)    \leq  C_\eta \frac{e^{6 \eta t} (C^{*}\, t^{\frac{1}{3}})}{\vert x \vert^{3}},$$ and then we have 
\begin{eqnarray}\label{eq78} \nonumber 
  I_a =  o(t)\, \left(\frac{1}{\vert x \vert^2}\right), \quad \text{when}\quad \vert x \vert \longrightarrow+\infty.
\end{eqnarray}
We study now the term $I_b$ in the formula (\ref{eq72}). By the estimate (\ref{eq77}) we get 
\begin{equation*}
I_b \leq \int_{0}^{t} \int_{\vert y \vert > \frac{\vert x \vert}{2}} \left\vert K_{\eta}(t-\tau, x-y)\right\vert  \left\vert u(\tau, y ) \partial_y u (\tau, y)\right\vert  dy\, d\tau \leq  \int_{0}^{t} \int_{\vert y \vert > \frac{\vert x \vert}{2}} \vert K_{\eta}(t-\tau, x-y) \vert \frac{C^{*}}{\tau^{\frac{2}{3}}(1+\vert y \vert^4)} dy d\tau, 
\end{equation*} but, since in the term $I_b$ above we have $\vert y \vert >\frac{\vert x \vert}{4}$ then  we can write 
$ \ds{\frac{1}{1+\vert y \vert^4}} \leq \frac{c}{\vert x \vert^4}$, and thus we obtain
\begin{eqnarray*}
 & & \int_{0}^{t} \int_{\vert y \vert > \frac{\vert x \vert}{2}} \vert K_{\eta}(t-\tau, x-y) \vert \frac{C^{*}}{\tau^{\frac{2}{3}}(1+\vert y \vert^4)} dy d\tau \leq \frac{C^{*}}{\vert x \vert^4} \int_{0}^{t} \int_{\vert y \vert >\frac{\vert x \vert}{4}} \vert K_{\eta}(t-\tau, x-y)\vert dy d \tau \\
 & \leq & \frac{C^{*}}{\vert x \vert^4} \int_{0}^{t} \Vert K_\eta(t-\tau, \cdot)\Vert_{L^1} d \tau, 
\end{eqnarray*} where,   by the estimate  (\ref{eq79}) we write 
\begin{equation*}
\frac{C^{*}}{\vert x \vert^4} \int_{0}^{t} \Vert K_\eta(t-\tau, \cdot)\Vert_{L^1} d \tau \leq \frac{C^{*}}{\vert x \vert^4}\int_{0}^{t} \left(c_\eta \frac{e^{5\eta (t-\tau)}}{(t-\tau)^{\frac{1}{3}}}\right) d\tau \leq \frac{C^{*}}{\vert x \vert^4} \left(c_\eta e^{5\eta t} t^{\frac{2}{3}} \right). 
\end{equation*} Then, for $\vert x \vert$ large enough we have $\ds{I_b \leq \frac{C^{*}}{\vert x \vert^4} \left(c_\eta e^{5\eta t} t^{\frac{2}{3}} \right)}$ and thus we can write 
\begin{equation}\label{eq80}
I_b = o(t)\, \left(\frac{1}{\vert x \vert^2}\right), \quad \text{when}\quad \vert x \vert \longrightarrow +\infty.
\end{equation}
We study the term $I_c$ in the equation (\ref{eq72}). By the estimates (\ref{eq64}) and (\ref{eq77})  we have 
\begin{eqnarray*}
I_c & \leq& \int_{0}^{t} \vert K_{\eta} (t-\tau, x)\vert  \left( \int_{\vert y \vert > \frac{\vert x \vert}{2}} \vert  u(\tau, y ) \partial_y u (\tau, y) \vert  dy \right) d\tau\\ &\leq&  \int_{0}^{t} \left(c_\eta \frac{e^{5\eta(t-\tau)}}{(t-\tau)^{\frac{1}{3}}} \frac{1}{1+\vert x \vert^2} \right) \left( \int_{\vert y\vert >\frac{\vert x \vert}{2} }\frac{C^{*}}{\tau^{\frac{2}{3}} (1+\vert y \vert^{4})}dy\right) \,d \tau\\
&\leq & \int_{0}^{t}\left(c_\eta  \frac{e^{5\eta t }}{(t-\tau)^{\frac{1}{3}}} \frac{1}{\vert x \vert^2}\right) \left(\int_{\vert y\vert >\frac{\vert x \vert}{2} } \frac{C^{*}}{\tau^{\frac{2}{3}}(1+\vert y \vert^2)(1+\vert y \vert^2) }dy \right)  \, d\tau=(a), 
\end{eqnarray*} but, remark that in the term $I_b$ above we have $\vert y \vert >\frac{\vert x \vert}{4}$, then  we can write 
$ \ds{\frac{1}{1+\vert y \vert^2}} \leq \frac{c}{\vert x \vert^2}$ and thus we get 
\begin{eqnarray*}
(a)	 & \leq &   \int_{0}^{t}\left(c_\eta  \frac{e^{5\eta t }}{(t-\tau)^{\frac{1}{3}}} \frac{1}{\vert x \vert^2}\right) \left(\int_{\vert y\vert >\frac{\vert x \vert}{2} } \frac{C^{*}}{\tau^{\frac{2}{3}}\vert x \vert^2(1+\vert y \vert^2) }dy \right)  \, d\tau \\
 & \leq &\frac{c_\eta e^{5\eta t } C^{*}}{\vert x\vert^4}  \int_{0}^{t} \left(\frac{1}{(t-\tau)^{\frac{1}{3}}}\right) \left( \int_{\R} \frac{d y}{\tau^{\frac{2}{3}}(1+\vert y \vert^2) }\right)d\tau  \leq  \frac{c_\eta e^{5\eta t } C^{*}}{\vert x\vert^4} \left(\int_{0}^{t} \frac{d\tau }{(t-\tau)^{\frac{1}{3}} \tau^{\frac{2}{3}}}\right) \leq   \frac{c_\eta e^{5\eta t } C^{*}}{\vert x\vert^4}.
\end{eqnarray*} Thus, for $\vert x \vert$ large enough we have the estimate $\ds{I_c \leq \frac{c_\eta e^{5\eta t } C^{*}}{\vert x\vert^4}}$ and then we can write 
\begin{equation}\label{eq81}
I_c = o(t)\, \left(\frac{1}{\vert x \vert^2}\right),\quad  \text{when} \quad \vert x \vert \longrightarrow +\infty.
\end{equation} Finally,  by the estimates given in formulas  (\ref{eq78}), (\ref{eq80}) and (\ref{eq81}), we can write the  estimate (\ref{eq73}) and the Theorem \ref{Th:asymptotics} is now proven. \finpv
%Corregido hasta aqui.
\subsection{Proof of Theorem \ref{Th:estim-below}}\label{sec:estim-below}
 For $t>0$ we write the solution $u(t,x)$ as the integral formulation (\ref{integral}) and we will start by the following estimates that we shall need later. For the linear term in (\ref{integral})  we write
\begin{eqnarray*}
K_\eta(t,\cdot)\ast u_0(x) &=&  \int_{\vert y \vert \leq \frac{\vert x \vert}{2}} (K_\eta(t,x-y)-K_\eta(t,x))u_0(y)dy+ K_\eta(t,x) \int_{\vert y \vert \leq \frac{\vert x \vert}{2}} u_0(y)dy \\
& & +  \int_{\vert y \vert > \frac{\vert x \vert}{2}} K_\eta(t,x-y)u_0(y)dy\\
& = & I_1+I_2+I_3,
\end{eqnarray*}	
where we will study the terms $I_1, I_2$ and $I_3$. Recall that the term $I_1$ was already treated in formulas (\ref{eq74}) and (\ref{eq106}) as follows: 
\begin{equation}\label{termI1}
I_1 \leq  \left( C_\eta \, \frac{e^{6\eta t}}{t^{\frac{2}{3}}} \Vert \,  \vert \cdot \vert u_0 \Vert_{L^1} \right) \frac{1}{1+\vert x \vert^3}.
\end{equation}
On the other hand, for the term $I_2$ we write
$$ I_2 = K_\eta(t,x) \left( \int_{\mathbb{R}}u_0(y)dy\right)   - K_\eta(t,x)  \left( \int_{\vert y \vert \geq \frac{\vert x \vert}{2}} u_0(y)dy\right),$$
where we will observe that the kernel $K_{\eta}(t,x)$ defined in (\ref{Kernel}) can be written for  $\vert x \vert$ large enough as $\ds{K_\eta(t,x)= \frac{c_{\eta,t}\,t}{\vert x \vert^2}}$, and where the quantity $c_{\eta,t}>0$ only depends on $\eta>0$ and $t>0$.  Indeed, by the identity (\ref{eq05}) and the identity (\ref{Descomposition-Kernel}) we can write  (for $\vert x \vert$ large)
$$ K_\eta(t,x)= \frac{2 \eta t}{(2 \pi i x)^2} +\frac{2 t}{(2 \pi i x)^2} (I_a + I_b) = \frac{1}{\vert x \vert^2} \times \frac{t}{-2 \pi^2} (  \eta  +  (I_a + I_b) ),$$ where: the quantity $I_a$ is given in expression (\ref{eq06}), the quantity $I_b$ is given in expression (\ref{eq07}), and moreover, by Lemma \ref{Lemma-tech1} we have $\ds{\vert I_a+I_b\vert \leq c_\eta  e^{4 \eta t}}$.  Thus, 
 we define the quantity $c_{\eta,t}$ as follows:
\begin{equation*}
c_{\eta,t} = \frac{1}{-2 \pi^2} (  \eta  +  (I_a + I_b)),
\end{equation*} and  we have the identity 
$$ K_\eta (t,x)= \frac{c_{\eta,t} \, t}{\vert x \vert^2}.$$ 
Once we dispose of this identity, the term $I_2$ is written as: 
\begin{equation}\label{termI2}
 I_2 = \frac{c_{\eta,t} \, t}{\vert x \vert^2} \left( \int_{\mathbb{R}}u_0(y)dy\right)   - \frac{c_{\eta,t} \, t}{\vert x \vert^2}  \left( \int_{\vert y \vert \geq \frac{\vert x \vert}{2}} u_0(y)dy\right).
\end{equation}  
Finally, remark that the term $I_3$ was already studied in formula (\ref{eq107}), and recalling this formula we have 
\begin{equation}\label{termI3}
 I_3 \leq \frac{c}{\vert x \vert^{2+\varepsilon}} \left( c_\eta \frac{e^{5 \eta t }}{t^{\frac{1}{3}}} \right).\\
\end{equation}
Now, in order to study the  nonlinear term in (\ref{integral}),  for $0<s<t$ we will start by studying the expression $\ds{ \partial_x K_\eta(t-s, \cdot) \ast u^{2}(s,x)}$.  Recall that  by point $2)$ of Lemma \ref{estim-K-2} we have  $$\ds{\vert \partial_x K_\eta(t,x)\vert  \leq C_\eta \frac{e^{6\eta t }}{t^{\frac{2}{3}}} \frac{1}{1+\vert x \vert^3}},$$ and moreover, by the estimate (\ref{decay}) given in Theorem \ref{Th-decay} we have 
$$ \vert u(s,x)\vert \leq \frac{C(s,\eta,u_0,u)}{1+\vert x \vert^2},$$
where the constant $C(s,\eta,u_0,u)>0$ (given in the formula (\ref{Constante-C})) is written as $\ds{ C(s,\eta,u_0,u)= \frac{\mathfrak{C}(T,\eta,u_0,u)}{s^{\frac{1}{3}}}}$.\\
\\
With these estimates in mind, we write now 
\begin{eqnarray*}
 \vert  \partial_x K_\eta(t-s, \cdot) \ast u^{2}(s,x) \vert & \leq & C_\eta \frac{e^{6\eta (t-s) }}{(t-s)^{\frac{2}{3}}} \frac{\mathfrak{C}^{2}(T,\eta,u_0,u)}{s^{\frac{2}{3}}}  \int_{\mathbb{R}} \frac{dy}{(1+\vert x-y\vert^{3})(1+ \vert y \vert^4)} \\
 &   \leq &  C_\eta \frac{e^{6\eta t }}{(t-s)^{\frac{2}{3}}} \frac{\mathfrak{C}^{2}(T,\eta,u_0,u)}{s^{\frac{2}{3}}} \frac{1}{1+\vert x \vert^3}.
\end{eqnarray*}
Hence, we obtain 
\begin{eqnarray} \label{nonlinterm} \nonumber
\left\vert \int_{0}^{t} K_\eta(t-s, \cdot) \ast \partial_x (u^{2}(s,x))ds \right\vert &\leq  &  C_\eta  e^{6\eta t} \mathfrak{C}^{2}(T,\eta,u_0,u) \left( \int_{0}^{t} \frac{ds}{(t-s)^{\frac{2}{3}} s^{\frac{2}{3}}}\right)  \frac{1}{1+\vert x \vert^3}\\
& \leq & C_\eta  e^{6\eta t} \frac{\mathfrak{C}^{2}(T,\eta,u_0,u) }{t^{\frac{1}{3}}} \frac{1}{1+\vert x \vert^3}.
\end{eqnarray}   
As mentioned before, these estimates given on the linear and the nonlinear term will be very useful to prove this theorem. We start by getting  back to the integral formulation (\ref{integral}) and we write the following profile for the solution: 
\begin{equation}\label{profile}
u(t,x)= I_1+ \underbrace{\frac{c_{\eta,t} \, t}{\vert x \vert^2} \left( \int_{\mathbb{R}}u_0(y)dy\right)  - \frac{c_{\eta,t} \, t}{\vert x \vert^2}  \left( \int_{\vert y \vert \geq \frac{\vert x \vert}{2}} u_0(y)dy\right)}_{=I_2}   + I_3 + \int_{0}^{t} \partial_x K_\eta(t-s, \cdot) \ast u^{2}(s,x)ds,
\end{equation} where we will consider the following cases:
\begin{enumerate}
\item[1)] The case $\ds{\int_{\mathbb{R}} u_0(y)dy\neq 0}$.  Remark  that once we dispose of the estimates for the terms $I_1, I_2$ and $I_3$, given in formulas (\ref{termI1}), (\ref{termI2}) and (\ref{termI3}) respectively,   we can write 
\begin{equation}\label{linterm}
I_1+ \frac{c_{\eta,t} \, t}{\vert x \vert^2} \left( \int_{\mathbb{R}}u_0(y)dy\right) - \frac{c_{\eta,t} \, t}{\vert x \vert^2}  \left( \int_{\vert y \vert \geq \frac{\vert x \vert}{2}} u_0(y)dy\right) + I_3 = \frac{c_{\eta,t} \, t}{\vert x \vert^2} \left( \int_{\mathbb{R}}u_0(y)dy\right)   +  o(t)\left( \frac{1}{\vert x \vert^2}\right), \quad \vert x \vert \longrightarrow +\infty.
\end{equation} On the other hand, for the nonlinear term, by  the estimate (\ref{nonlinterm}) we can write 
\begin{equation*}
\int_{0}^{t} \partial_x K_\eta(t-s,\cdot)\ast  u^2(s,x) ds = o(t) \left(\frac{1}{\vert x \vert^2}\right), \quad \vert x \vert \longrightarrow +\infty. 
\end{equation*}
Thus, getting back to the profile  (\ref{profile}), for $\vert x \vert$ large enough we can write:
\begin{eqnarray*}
	\vert u(t,x) \vert &=& \left\vert  \frac{c_{\eta,t} \, t}{\vert x \vert^2} \left( \int_{\mathbb{R}}u_0(y)dy\right)   +  o(t)\left( \frac{1}{\vert x \vert^2}\right) \right\vert=  \left\vert  \frac{c_{\eta,t} \, t}{\vert x \vert^2} \left( \int_{\mathbb{R}}u_0(y)dy\right)  - \left( -  o(t)\left( \frac{1}{\vert x \vert^2}\right) \right) \right\vert \\
	& \geq & \frac{c_{\eta,t} \, t}{\vert x \vert^2} \left\vert \int_{\mathbb{R}}u_0(y)dy \right\vert - \left\vert o(t)\left( \frac{1}{\vert x \vert^2}\right) \right\vert. 
\end{eqnarray*}
Then, recalling the definition of the quantity $ \ds{ o(t)\left(\frac{1}{\vert x \vert^2}\right)}$ given in the formula (\ref{little-o}), for $\ds{ \frac{c_{\eta,t} \, t}{2}   \left\vert \int_{\mathbb{R}}u_0(y)dy \right\vert>0}$ there exists $M>0$ such that for all $\vert x \vert >M$ we have $\ds{\left\vert     o(t)\left(\frac{1}{\vert x \vert^2}\right) \right\vert \leq \frac{c_{\eta,t} \, t}{2 \vert x \vert^2}  \left\vert \int_{\mathbb{R}}u_0(y)dy \right\vert}$. Hence, we have
$$ - \left\vert     o(t)\left(\frac{1}{\vert x \vert^2}\right) \right\vert \geq -  \frac{c_{\eta,t} \, t}{2 \vert x \vert^2}  \left\vert \int_{\mathbb{R}}u_0(y)dy \right\vert,$$ and getting back to the estimate from below on the quantity $\vert u (t,x)\vert$ above we obtain 
$$ \vert u(t,x)\vert \geq  \frac{c_{\eta,t} \, t}{2 \vert x \vert^2}  \left\vert \int_{\mathbb{R}}u_0(y)dy \right\vert.$$
\item[2)] The case $\ds{\int_{\mathbb{R}} u_0(y)dy=0}$. Remark  that, always by the estimates given on the terms $I_1, I_2$ and $I_3$ (see (\ref{termI1}), (\ref{termI2}) and (\ref{termI3}))  we can write now
\begin{eqnarray*}
& &  I_1+ \underbrace{\frac{c_{\eta,t} \, t}{\vert x \vert^2} \left( \int_{\mathbb{R}}u_0(y)dy\right)   - \frac{c_{\eta,t} \, t}{\vert x \vert^2}  \left( \int_{\vert y \vert \geq \frac{\vert x \vert}{2}} u_0(y)dy\right)}_{=I_2}  + I_3 \\
 & \leq  &  \frac{c_{\eta,t} \, t}{\vert x \vert^2}  \left( \int_{\vert y \vert \geq \frac{\vert x \vert}{2}} \vert u_0(y) \vert dy\right) +  \frac{c}{\vert x \vert^{2+\varepsilon}} \left( c_\eta \frac{e^{5 \eta t }}{t^{\frac{1}{3}}} \right) \\ \nonumber  
& \leq & \frac{c_{\eta,t} \, t}{\vert x \vert^{2+\varepsilon}}  \left( \int_{\vert y \vert \geq \frac{\vert x \vert}{2}} \vert x \vert^{\varepsilon} \vert u_0(y) \vert dy\right) + \frac{c}{\vert x \vert^{2+\varepsilon}} \left( c_\eta \frac{e^{5 \eta t }}{t^{\frac{1}{3}}} \right)  \\  \nonumber  
&\leq &  \frac{c_{\eta,t} \, t}{\vert x \vert^{2+\varepsilon}}  \left( \int_{\vert y \vert \geq \frac{\vert x \vert}{2}} \vert y \vert^{\varepsilon} \vert u_0(y) \vert dy\right) + \frac{c}{\vert x \vert^{2+\varepsilon}} \left( c_\eta \frac{e^{5 \eta t }}{t^{\frac{1}{3}}} \right) \\
&\leq & \frac{c}{\vert x \vert^{2+\varepsilon}} \left(\Vert \, \vert \cdot \vert^{\varepsilon}\, u_0 \Vert_{L^1}+ c_\eta \frac{e^{5 \eta t }}{t^{\frac{1}{3}}}\right).
\end{eqnarray*}  
With this estimate and the estimate for the nonlinear term given in the formula (\ref{nonlinterm})  we can write 
\begin{equation*}
\vert u(t,x)\vert \leq \frac{c}{\vert x \vert^{2+\varepsilon}} \left(\Vert \, \vert \cdot \vert^{\varepsilon}\, u_0 \Vert_{L^1}+ c_\eta \frac{e^{5 \eta t }}{t^{\frac{1}{3}}}\right)+C_\eta  e^{6\eta t} \frac{\mathfrak{C}^{2}(T,\eta,u_0,u) }{t^{\frac{1}{3}}} \frac{1}{1+\vert x \vert^3}, 
\end{equation*} and for $\vert x \vert $ large enough and $0<\varepsilon \leq 1$ we have 
\begin{equation*}
\vert u(t,x)\vert \leq \frac{c}{1+\vert x \vert^{2+\varepsilon}} \left( \Vert \, \vert \cdot \vert^{\varepsilon}\, u_0 \Vert_{L^1}+ \frac{c_\eta e^{5 \eta t}  + C_\eta  e^{6\eta t}\mathfrak{C}^{2}(T,\eta,u_0,u)}{t^{\frac{1}{3}}}\right).
\end{equation*}
Theorem \ref{Th:estim-below} is proven.  \finpv 
\end{enumerate} 
%%%%%%%%%%%%%%%%%%%%%%%%%%%%%%%%%%%%%%%%%%%
%%%%%%%%%%%%%%%%%%%%%%%%%%%%%%%%%%%%%%%%%%%
\section{The LWP in Lebesgue spaces: proof of Theorem \ref{Th-Lebesgue}}\label{sec:Lebesgue} 
We start by remarking that the kernel $K_\eta(t,\cdot)$ given in (\ref{Kernel}) and its derivative  $\partial_x K_\eta(t,\cdot)$ belong to the space $L^{p}(\R)$ for $1\leq p \leq +\infty$. Indeed,  by point $1)$ of Proposition \ref{Prop1} we have, for all time $t>0$, $\ds{\vert K_\eta(t,x)\vert \leq c_\eta  \frac{e^{5 \eta t }}{t^{\frac{1}{3}}} \frac{1}{1+\vert x \vert^2}}$ and then, for $1\leq p \leq +\infty$ we get $\ds{\Vert K_{\eta}(t,\cdot) \Vert_{L^p} \leq  c_\eta   \frac{e^{5 \eta t }}{t^{\frac{1}{3}}} \left\Vert \frac{1}{1+\vert \cdot \vert^2} \right\Vert_{L^p}}$, hence, for the sake of simplicity, we will write 
\begin{equation}\label{eq100}
\Vert K_\eta(t,\cdot)\Vert_{L^p} \leq  c_\eta   \frac{e^{5 \eta t }}{t^{\frac{1}{3}}}. 
\end{equation}
In the same way, recall that by point $2)$ of Lemma \ref{estim-K-2}, we have, for all time $t>0$, $\ds{\vert \partial_x K_\eta(t,x)\vert \leq C_\eta \frac{e^{6\eta t}}{t^{\frac{2}{3}}} \frac{1}{1+\vert x \vert^3}}$, thereby, for $1\leq p \leq +\infty$,   we obtain  
\begin{equation}\label{eq101}
\Vert \partial_x K_\eta(t,\cdot) \Vert_{L^p} \leq C_\eta \frac{e^{6 \eta t}}{t^{\frac{2}{3}}}. 
\end{equation} 
Estimates (\ref{eq100}) and (\ref{eq101}) will allow us to study the existence of \emph{mild} solutions for the Cauchy problem (\ref{eq:f}) in the framework of Lebesgue spaces and when the initial datum $u_0$ is small enough.  Let $T>0$ fix and consider  the Banach space $L^{\infty}(]0, T[, L^{p}(\R))$ with the norm  $\ds{\sup_{0<t<T} t^{\frac{1}{3}} \Vert \cdot \Vert_{L^{p}_{x}}}$. We write 
\begin{equation*}
\sup_{0<t<T}t^{\frac{1}{3}} \Vert u(t,\cdot) \Vert_{L^p} \leq \sup_{0<t<T}t^{\frac{1}{3}} \Vert  K_\eta(t,\cdot) \ast u_0 \Vert_{L^p} + \sup_{0<t<T}t^{\frac{1}{3}} \left\Vert  \int_{0}^{t} K_\eta(t-s,\cdot)\ast \partial_{x} (u^{2}(s,\cdot)) ds \right\Vert_{L^p},
\end{equation*} and we will estimate each term nn the right-hand side. \\
\\
For the first term, by the  estimate (\ref{eq100}) we can write 
\begin{equation}\label{eq105}
\sup_{0<t<T}t^{\frac{1}{3}} \Vert  K_\eta(t,\cdot) \ast u_0 \Vert_{L^p} \leq \sup_{0<t<T}t^{\frac{1}{3}} \Vert  K_\eta(t,\cdot) \Vert_{L^1} \Vert u_0 \Vert_{L^p} \leq  \sup_{0<t<T}t^{\frac{1}{3}}  \left( c_\eta   \frac{e^{5 \eta t }}{t^{\frac{1}{3}}} \right) \Vert u_0 \Vert_{L^p} \leq c_\eta e^{5 \eta T} \Vert u_0 \Vert_{L^p}. 
\end{equation} 
Now,  the second term  is estimated as follows: first, for all time  $t\in ]0,T[$ and for $1\leq q \leq +\infty$ which verifies $1+\frac{1}{p}= \frac{1}{q}+ \frac{2}{p}$, we write 
\begin{eqnarray*}
	\left\Vert  \int_{0}^{t} K_\eta(t-s)\ast \partial_{x} (u^{2}(s,\cdot)) ds \right\Vert_{L^p} & \leq& \int_{0}^{t} \Vert K_\eta(t-s)\ast \partial_{x} (u^{2}(s,\cdot)) \Vert_{L^p} ds \leq \int_{0}^{t}  \Vert \partial_x K_{\eta}(t-s,\cdot)\ast u^{2}(s,\cdot) \Vert_{L^p} ds\\
	&\leq &   \int_{0}^{t}  \Vert \partial_x K_{\eta}(t-s,\cdot)\Vert_{L^q} \Vert  u^{2}(s,\cdot) \Vert_{L^\frac{p}{2}} ds, 
\end{eqnarray*} and then, by the estimate (\ref{eq101}) we get 
\begin{eqnarray*}
	\int_{0}^{t}  \Vert \partial_x K_{\eta}(t-s,\cdot)\Vert_{L^q} \Vert  u^{2}(s,\cdot) \Vert_{L^\frac{p}{2}} ds &\leq & \int_{0}^{t}  \left( C_\eta \frac{e^{6\eta (t-s)}}{(t-s)^{\frac{2}{3}}} \right)   \Vert  u^{2}(s,\cdot) \Vert_{L^\frac{p}{2}} ds \leq C_{\eta} e^{6\eta T} \int_{0}^{t} \frac{1}{(t-s)^{\frac{2}{3}}} \Vert u(s,\cdot)\Vert^{2}_{L^p}ds \\
	&\leq & C_{\eta} e^{6\eta T} \int_{0}^{t} (t-s)^{-\frac{2}{3}} s^{-\frac{2}{3}} \left( s^{\frac{1}{3}}  \Vert u(s,\cdot)\Vert_{L^p}\right)^{2} ds \\
	&\leq & C_{\eta} e^{6\eta T} \left( \sup_{0<t<T} t^{\frac{1}{3}} \Vert u(t,\cdot)\Vert_{L^p} \right)^{2} \left(  \int_{0}^{t} (t-s)^{-\frac{2}{3}} s^{-\frac{2}{3}}ds\right), 
\end{eqnarray*} where, the last expression (also known as the Beta function) verifies 
$\ds{ \int_{0}^{t} (t-s)^{-\frac{2}{3}} s^{-\frac{2}{3}}ds \leq c t^{-\frac{1}{3}}}$, and thus we can write 
\begin{equation*}
\left\Vert  \int_{0}^{t} K_\eta(t-s)\ast \partial_{x} (u^{2}(s,\cdot)) ds \right\Vert_{L^p} \leq C_{\eta} e^{6\eta T} \left( \sup_{0<t<T_0} t^{\frac{1}{3}} \Vert u(t,\cdot)\Vert_{L^p} \right)^{2} t^{-\frac{1}{3}}. 
\end{equation*} Once we have this estimate we write 
\begin{eqnarray}\label{eq104} \nonumber 
\sup_{0<t<T} t^{\frac{1}{3}}  \left\Vert  \int_{0}^{t} K_\eta(t-s)\ast \partial_{x} (u^{2}(s,\cdot)) ds \right\Vert_{L^p} &\leq& \sup_{0<t<T} t^{\frac{1}{3}} \left( C_{\eta} e^{6\eta T} \left( \sup_{0<t<T_0} t^{\frac{1}{3}} \Vert u(t,\cdot)\Vert_{L^p} \right)^{2} t^{-\frac{1}{3}} \right) \\
&\leq & C_\eta e^{6\eta T} \left( \sup_{0<t<T} t^{\frac{1}{3}} \Vert u(t,\cdot)\Vert_{L^p} \right)^{2}. 
\end{eqnarray}
Now, with the estimates (\ref{eq105}) and (\ref{eq104}) we set the quantity $\delta$ as $\ds{\delta = \frac{1}{4 c_\eta C_\eta e^{11 \eta T}}>0}$, and if  the initial datum verifies $\Vert u_0 \Vert_{L^p} < \delta$ then we apply the  Picard contraction principle to obtain  a solution $u(t,x)$ of the integral equation (\ref{integral}).  \\
\\
We prove now the uniqueness of this solution  and for this we will follow the same ideas given at the end of the proof of Theorem \ref{Th-aux-1}. Indeed, let us suppose that the integral equation  (\ref{integral}) admits two solutions  $u_1$  and $u_2$ (arising  from the same initial datum $u_0$) in the space $\left(L^{\infty}(]0, T[, L^{p}(\R)), \ds{\sup_{0<t<T} t^{\frac{1}{3}} \Vert \cdot \Vert_{L^{p}_{x}}}\right)$. Then, we denote $v=u_1-u_2$ and by recalling the identity (\ref{v}) we write 
$$ v(t,\cdot)=-\frac{1}{2} \int_{0}^{t} K_\eta(t-s, \cdot)\ast \left(  \partial_{x} (v(s,\cdot) u_1(s,\cdot)+u_2(s,\cdot) v(s,\cdot)) \right)ds. $$
Finally, let $0<T^{*}\leq T$ be the maximal time such that we have $v=0$ on the interval  $]0,T^{*}[$ and we will prove that $T^{*}=T$. Indeed, if we suppose   $T^{*}<T$ then there exists a time $T^{*}<T_1<T$ and we consider the space $L^{\infty}(]T^{*}, T_1[,L^{p}(\R))$ endowed with the norm $\Vert \cdot \Vert_{(T_1-T^{*})}= \ds{\sup_{T^{*}<t<T_1} t^{\frac{1}{3}} \Vert \cdot \Vert_{L^{p}_{x}}}$. Now, remark that by the estimate (\ref{eq104}) we can write 
$$ \Vert v \Vert_{(T_1-T^{*})} \leq  C_\eta e^{6\eta (T_1- T^{*}) } \Vert v \Vert_{(T_1-T^{*})} \left( \Vert u_1 \Vert_{(T_1-T^{*})} + \Vert u_2 \Vert_{(T_1-T^{*})} \right),$$ and remark also that for a function $f \in L^{\infty}(]T^{*}, T_1[,L^{p}(\R))$ we have $\ds{\lim_{T_1 \longrightarrow T^{*}} \Vert f \Vert_{(T_1-T^{*})}=0}$. Therefore, we can take $0<T_1-T^{*}$ small enough such that 
$$ \Vert u_1 \Vert_{(T_1-T^{*})} + \Vert u_2 \Vert_{(T_1-T^{*})} \leq \frac{1}{2\, C_\eta e^{6\eta (T_1- T^{*}) }}.$$ By this inequality and the previous estimate on the quantity $\Vert v \Vert_{(T_1-T^{*})}$ we obtain $ \Vert v \Vert_{(T_1-T^{*})} =0$ which contradicts the definition of $T^{*}$.  Theorem  \ref{Th-Lebesgue} is now proven. \finpv
\section{Appendix}
\subsection*{Proof of Lemma \ref{Lemma-tech1}}
Recall that the  term $I_a$  in \eqref{eq06} is  given as  
\begin{eqnarray*}
I_a &=& \int_{\xi <0} e^{2\pi i x \xi} \partial_{\xi} \left( (e^{i t \xi^3- \eta t (-\xi^3 +\xi)})(3i t \xi^2 -\eta t (-3\xi^2+1)) \right) d \xi = \int_{\xi <0} e^{2\pi i x \xi} \partial_{\xi} \left( \partial_{\xi} (e^{i t \xi^3- \eta t (-\xi^3 +\xi)})  \right) d \xi \\
&=& \int_{\xi <0}  e^{2\pi i x \xi}  \ \partial^{2}_{\xi} \left( e^{i t \xi^3- \eta t (-\xi^3 +\xi)} \right) d \xi = \int_{\xi <0}  e^{2\pi i x \xi}  \ \partial^{2}_{\xi} \widehat{K}_\eta (t , \xi)  \  d \xi.
\end{eqnarray*}
On the other hand,  by Lemma 5.1 in  \cite{BorysAlvarez-tesis}, we have  for all $\xi \neq 0$:
$$  \partial_{\xi}^{2}\widehat{K}_\eta (t , \xi) = \widehat{K}_\eta (t , \xi) t^{2} \Big( 3i\xi^{2}  - \eta  \ \sign(\xi) (3 \xi^{2}- 1) \Big)^{2}+ 6t \xi (i - \eta  \  \sign(\xi)) \widehat{K}_\eta (t , \xi),$$ and then we can write
\begin{equation} \label{Ap1}
\begin{split}
| I_a | \leq  \left \| \partial^2_{\xi}  \widehat{K}_\eta (t , \cdot ) \right\|_{L^{1}(]-\infty, 0[)}  \leq c(1+\eta)^2 \  t^2  \left \| \widehat{K}_\eta (t , \cdot )  (1 + \left | \cdot \right|^4) \right \|_{L^{1}(\mathbb{R})} \\
 + c(1+\eta)  \ t  \left \| \widehat{K}_\eta (t , \cdot ) (1 + \left | \cdot \right|) \right \|_{L^{1}(\mathbb{R})}.
\end{split}
\end{equation}
In order to study the term on the right-hand side we have the following estimates:  for $m> -1$, by the estimate (\ref{Estim-K-Fou-L1}) and denoting by $\Gamma$  the ordinary gamma function, we have:
\begin{eqnarray}\label{Estim-m} \nonumber
\left \|(1 +\left | \cdot\right |^{m}) \widehat{K}_\eta (t , \cdot ) \right \|_{L^{1}} &\leq&   \left \| \widehat{K}_\eta (t , \cdot ) \right \|_{L^{1}} +  \left \| \left | \xi\right |^{m} \widehat{K}_\eta (t , \cdot ) \right \|_{L^{1}} \\ \nonumber 
&\leq& C \frac{e^{3\eta t }}{(\eta t)^{\frac{1}{3}}}   +  \displaystyle\int_{\left |\xi \right | \leq 2}  \left | \xi \right |^{m} e^{-t\eta(\left |\xi\right|^{3}-\left |\xi\| \right)} \, d\xi +  \displaystyle\int_{\left |\xi \right | \geq 2}  \left | \xi \right |^{m}  e^{-t\eta \frac{3}{4}\left |\xi\right|^{3}} \, d\xi  \\ \nonumber 
&\leq&   C \frac{e^{3\eta t }}{(\eta t)^{\frac{1}{3}}} + \frac{2^{m+2}}{m+1}  \ e^{2\eta t } +  \frac{c_m \Gamma(\frac{m+1}{3})}{(\eta t )^{(\frac{m+1}{3})}}  \\
&\leq&  C_{m} \frac{e^{3\eta t }}{(\eta t)^{\frac{1}{3}}} + C_m \frac{1}{(\eta t)^\frac{m+1}{3}}. 
\end{eqnarray}
With this estimate  (setting first $m=4$ and then $m=1$) we get back to \eqref{Ap1} and we write 
\begin{eqnarray}\label{eq92} \nonumber 
\vert I_a \vert & \leq & c(1+\eta)^2 t^2 \left( \frac{e^{3\eta t }}{(\eta t)^{\frac{1}{3}}} +  \frac{1}{(\eta t)^\frac{5}{3}} \right)+ c(1+\eta) t \left( \frac{e^{3\eta t }}{(\eta t)^{\frac{1}{3}}} + \frac{1}{(\eta t)^\frac{2}{3}}\right) \\ \nonumber
&\leq & c\frac{(1+\eta)^2}{\eta^2}  (\eta t)^2 \left( \frac{e^{3\eta t }}{(\eta t)^{\frac{1}{3}}} +  \frac{1}{(\eta t)^\frac{5}{3}} \right)+ c\frac{(1+\eta)}{\eta} (\eta t)  \left( \frac{e^{3\eta t }}{(\eta t)^{\frac{1}{3}}} + \frac{1}{(\eta t)^\frac{2}{3}}\right) \\ \nonumber
&\leq& c \frac{(1+\eta)^2}{\eta^2}  \left( (\eta t)^{\frac{5}{3}} e^{3 \eta t} + (\eta t)^{\frac{1}{3}} \right) + c\frac{(1+\eta)}{\eta} \left( (\eta t)^{\frac{2}{2}} e^{3 \eta t} + (\eta t)^{\frac{1}{3}}\right) \\ \nonumber
& \leq &  c \frac{(1+\eta)^2}{\eta^2}  \left(2  e^{4 \eta t}  \right) + c\frac{(1+\eta)}{\eta} \left( 2 e^{4 \eta t}  \right) \\ \nonumber 
&\leq &  c \left( \frac{1+\eta}{\eta}\right)\left( \left( \frac{1+\eta}{\eta}\right) +1 \right)e^{4\eta t}\\ \nonumber 
& \leq & c  \left( \left( \frac{1+\eta}{\eta}\right)+ 1 \right)\left( \left( \frac{1+\eta}{\eta}\right) +1 \right)e^{4\eta t} \\
& \leq & c \left(\frac{1}{\eta} + 2 \right)^2 e^{4\eta t}.
\end{eqnarray}
The term $I_b$  in \eqref{eq06}  is treated following the same computations done for the term $I_a$ above.  \finpv
\subsection*{Proof of Lemma \ref{estim-K-2}} 
\begin{enumerate}
\item[1)] Remark first that as we have $K_{\eta}(t,x)= \mathcal{F}^{-1}\left( e^{(i\xi^3t - \eta t(\vert \xi \vert^3-\vert \xi \vert))} \right)(x)$ and by the identity $\partial_x K_{\eta}(t,x)= \mathcal{F}^{-1}\left( (2\pi i \xi) e^{(i\xi^3t - \eta t(\vert \xi \vert^3-\vert \xi \vert))} \right)(x)$, and moreover, as  the function $\ds{e^{(i\xi^3t - \eta t(\vert \xi \vert^3-\vert \xi \vert))}}$ and the function $\ds{ (2\pi i \xi) e^{(i\xi^3t - \eta t(\vert \xi \vert^3-\vert \xi \vert))}}$ belong to the space $L^{1}(\R)$, then by the properties of the inverse Fourier transform we have that $\ds{K_{\eta}(t,x)}$ and $\ds{\partial_x K_{\eta}(t,x)}$ are continuous functions and thus we have $K_{\eta}(t,\cdot) \in \mathcal{C}^{1}(\R)$. \\
\\
Now, we write 
\begin{eqnarray*} 
\partial_x K_\eta(t,x) &=& \int_{\mathbb{R}} (2 \pi i \xi) e^{2\pi i x \xi} \widehat{K_{\eta}}(t,\xi) d \xi = \frac{1}{2 \pi i x } \int_{\xi <0} (2 \pi i \xi) (2\pi i x) e^{2 \pi i x \xi }\widehat{K_{\eta}}(t,\xi) d \xi \\
& &  + \frac{1}{2\pi i x}  \int_{\xi >0} (2 \pi i \xi) (2\pi i x) e^{2 \pi i x \xi }\widehat{K_{\eta}}(t,\xi) d \xi,
\end{eqnarray*} and since $\partial_{\xi}(e^{2\pi i x \xi})= 2\pi i x e^{2\pi i x \xi}$ then we can write 
\begin{eqnarray*}
& & \frac{1}{2 \pi i x } \int_{\xi <0} (2 \pi i \xi) (2\pi i x) e^{2 \pi i x \xi }\widehat{K_{\eta}}(t,\xi)d \xi+ \frac{1}{2\pi i x}  \int_{\xi >0} (2 \pi i \xi) (2\pi i x) e^{2 \pi i x \xi }\widehat{K_{\eta}}(t,\xi) d \xi \\
&=& \frac{1}{2 \pi i x } \int_{\xi <0} \partial_{\xi}(e^{2\pi i x \xi}) (2\pi i \xi ) e^{i t \xi^3- \eta t (-\xi^3 +\xi)} d \xi + \frac{1}{2\pi i x}  \int_{\xi >0} \partial_{\xi}(e^{2\pi i x \xi})(2\pi i \xi )e^{i t \xi^3- \eta t (\xi^3 -\xi)} d \xi, 
\end{eqnarray*}thereafter, we integrate by parts and  we get 
\begin{eqnarray}\label{eq98} \nonumber
& & \frac{1}{2 \pi i x } \int_{\xi <0} \partial_{\xi}(e^{2\pi i x \xi}) (2\pi i \xi ) \widehat{K_{\eta}}(t,\xi)d \xi + \frac{1}{2\pi i x}  \int_{\xi >0} \partial_{\xi}(e^{2\pi i x \xi})(2\pi i \xi )\widehat{K_{\eta}}(t,\xi)d \xi \\ \nonumber
&=&  \frac{1}{2 \pi i x } \int_{\xi <0}  e^{2\pi i x \xi}  (2\pi i) \widehat{K_{\eta}}(t,\xi) d \xi + \frac{1}{2\pi i x}  \int_{\xi >0} (e^{2\pi i x \xi})(2\pi i )\widehat{K_{\eta}}(t,\xi)d \xi \\ \nonumber
& & +  \frac{1}{2 \pi i x } \int_{\xi <0} e^{2\pi i x \xi}  (2\pi i \xi ) \partial_{\xi} \widehat{K_{\eta}}(t,\xi)d \xi + \frac{1}{2\pi i x}  \int_{\xi >0} e^{2\pi i x \xi}(2\pi i \xi )\partial_{\xi}\widehat{K_{\eta}}(t,\xi)d \xi \\ \nonumber 
&=& \frac{1}{x} \left( \int_{\xi <0}  e^{2\pi i x \xi}  \widehat{K_{\eta}}(t,\xi)d \xi +   \int_{\xi >0} (e^{2\pi i x \xi})\widehat{K_{\eta}}(t,\xi)d \xi  \right) \\ \nonumber
& & + \frac{1}{x} \left( \int_{\xi <0} e^{2\pi i x \xi}  \xi  \partial_{\xi}\widehat{K_{\eta}}(t,\xi) d \xi +   \int_{\xi >0} e^{2\pi i x \xi}  \xi \partial_{\xi}\widehat{K_{\eta}}(t,\xi) d \xi \right)  \\
&=& I_1 + I_2. 
\end{eqnarray}  
In order to study the term $I_1$ remark that  we have $\ds{I_1= \frac{1}{x} K_\eta(t,x)}$ and by the estimate (\ref{eq10}) we obtain
\begin{equation}\label{eq97}
\vert I_1 \vert  \leq C_\eta  \frac{e^{5 \eta t }}{\vert x \vert^3}.
\end{equation} We study now the term $I_2$ above. Remark that the have $\ds{\partial^{2}_{\xi}(e^{2\pi i x \xi}) = -4\pi^2 x^2 e^{2\pi i x \xi}}$, and therefore we write 
\begin{eqnarray*}
I_2&=& \frac{1}{(-4 \pi^2 x^2)x} \left( \int_{\xi <0} (-4 \pi x^2) e^{2\pi i x \xi}  \xi  \partial_{\xi} \widehat{K_{\eta}}(t,\xi) d \xi +   \int_{\xi >0} (-4 \pi x^2) e^{2\pi i x \xi}  \xi \partial_{\xi}\widehat{K_{\eta}}(t,\xi)d \xi \right)\\
&=& \frac{1}{-4 \pi^2 x^3} \left( \int_{\xi <0} \partial^{2}_{\xi} (e^{2\pi i x \xi})  \xi  \partial_{\xi} \widehat{K_{\eta}}(t,\xi) d \xi +   \int_{\xi >0} \partial^{2}_{\xi} (e^{2\pi i x \xi})  \xi \partial_{\xi}\widehat{K_{\eta}}(t,\xi)d \xi \right),
\end{eqnarray*}then, integrating   by parts the last expression  we can write 
\begin{equation}\label{eq95}
I_2 =\frac{1}{-4 \pi^2 x^3} \left(\underbrace{ \int_{\xi <0} e^{2\pi i x \xi} \left(2 \partial^{2}_{\xi} \widehat{K_{\eta}}(t,\xi) + \xi \partial^{3}_{\xi} \widehat{K_{\eta}}(t,\xi) \right)  d \xi}_{=(I_2)_a} +  \underbrace{\int_{\xi >0}  e^{2\pi i x \xi} \left(2 \partial^{2}_{\xi} \widehat{K_{\eta}}(t,\xi) + \xi \partial^{3}_{\xi} \widehat{K_{\eta}}(t,\xi) \right)   d \xi }_{=(I_2)_b} \right),
\end{equation} and now we will prove the following estimate
\begin{equation}\label{eq94}
\vert (I_2)_a \vert + \vert (I_2)_b \vert \leq C_\eta e^{5 \eta t}.
\end{equation}
%Revisado hasta aqui.
Indeed, for the term $\ds{(I_2)_a}$ we write $\ds{\vert (I_2)_a \vert \leq c \Vert \partial^{2}_{\xi} \widehat{K_{\eta}}(t,\cdot) \Vert_{L^{1}(]-\infty, 0 [)} + c\Vert \xi \partial^{3}_{\xi} \widehat{K_{\eta}}(t,\cdot) \Vert_{L^{1}(]-\infty, 0 [)} }$, but recall that by the  estimates  (\ref{Ap1}) and (\ref{eq92}) we have $\ds{\Vert \widehat{K_{\eta}}(t,\cdot) \Vert_{L^{1}(]-\infty, 0 [)} \leq C_\eta e^{4\eta t}}$ and therefore we can write 
\begin{equation}\label{eq93}
\vert (I_2)_a \vert \leq C_\eta e^{4 \eta t } + c\Vert \xi \partial^{3}_{\xi} \widehat{K_{\eta}}(t,\cdot) \Vert_{L^{1}(]-\infty, 0 [)} \leq  C_\eta e^{5 \eta t } + c\Vert \xi \partial^{3}_{\xi} \widehat{K_{\eta}}(t,\cdot) \Vert_{L^{1}(]-\infty, 0 [)}
\end{equation} Now, we study the term $\ds{c\Vert \xi \partial^{3}_{\xi}. \widehat{K_{\eta}}(t,\cdot) \Vert_{L^{1}(]-\infty, 0 [)}}$.  By Lemma 5.1 in  \cite{BorysAlvarez-tesis} we have for all $\xi \neq 0$:
\begin{eqnarray*}
\partial^{3}_{\eta}\widehat{K_{\eta}}(t,\xi)&=& t^{3} \widehat{K_{\eta}}(t,\xi)(3 i \xi^2-\eta \sign(\xi)(3 \xi^2 -1))^{3} \\
& &   +t^2 \widehat{K_{\eta}}(t,\xi)(36 \xi^3(\eta^2-1)-72\,i \,\eta\, \sign(\xi) \xi^3+12\, i\, \eta \,\sign(\xi) \xi -12 \eta^2 \xi) \\
& &+ 6 t^{2}\widehat{K_{\eta}}(t,\xi) (\xi (i-\eta\, \sign(\xi)))(3i \xi^2 -\eta \,\sign(\xi)(3 \xi^2-1))  + 6t \widehat{K_{\eta}}(t,\xi) (i-\eta \, \sign(\xi)).   
\end{eqnarray*} Thus, we can write 
\begin{equation*}
\vert \partial^{3}_{\eta}\widehat{K_{\eta}}(t,\xi) \vert  \leq  C_\eta t^3 (1+\vert \xi \vert^6) \vert \widehat{K_{\eta}}(t,\xi)\vert + C_\eta t^2 (1+\vert \xi \vert^3) \vert \widehat{K_{\eta}}(t,\xi)\vert	+ C_\eta t \vert \widehat{K_{\eta}}(t,\xi) \vert,  
\end{equation*}  and thus we get 
\begin{equation*}
\vert \xi \vert \vert \partial^{3}_{\eta}\widehat{K_{\eta}}(t,\xi) \vert  \leq  C_\eta t^3 (1+\vert \xi \vert^7) \vert \widehat{K_{\eta}}(t,\xi)\vert + C_\eta t^2 (1+\vert \xi \vert^4) \vert \widehat{K_{\eta}}(t,\xi)\vert	+ C_\eta t (1+\vert \xi \vert) \vert \widehat{K_{\eta}}(t,\xi) \vert. 
\end{equation*} With this estimate we can write 

\begin{eqnarray*}
\Vert \xi \partial^{3}_{\xi} \widehat{K_{\eta}}(t,\cdot) \Vert_{L^{1}(]-\infty, 0[)} & \leq & \Vert \xi \partial^{3}_{\xi} \widehat{K_{\eta}}(t,\cdot) \Vert_{L^{1}(\R)} \leq c_\eta t^{3} \Vert (1+\vert \xi \vert^7) \widehat{K_{\eta}}(t,\cdot)\Vert_{L^{1}(\R)}  \\
& &   + c_\eta t^2 \Vert (1+\vert \xi \vert^4) \widehat{K_{\eta}}(t,\cdot) \Vert_{L^{1}(\R)}  + C_\eta t  \Vert (1+\vert \xi \vert) \widehat{K_{\eta}}(t,\cdot) \Vert_{L^{1}(\R)} \\
&=& (a), 
\end{eqnarray*}  and then, by the estimate (\ref{Estim-m})  (setting first $m=7$, thereafter  $m=4$ and finally $m=1$) we have  

\begin{eqnarray*}
(a)& \leq & C_\eta t^{3} \left( e^{2t\eta} + t^{-1/3} + t^{-(\frac{8}{3})} \right) + c_\eta t^{2} \left(e^{2t\eta} + t^{-1/3} + t^{-(\frac{5}{3})} \right)  + c_\eta t \left( e^{2t\eta} + t^{-1/3} + t^{-(\frac{2}{3})} \right)\\
&\leq & C_\eta e^{5 \eta t},
\end{eqnarray*} and thus we can write $\ds{\Vert \xi \partial^{3}_{\xi} \widehat{K_{\eta}}(t,\cdot) \Vert_{L^{1}(]-\infty, 0[)} \leq C_\eta e^{5\eta t}}$. With this estimate,  we get back to the estimate (\ref{eq93}) and we write $\ds{\vert (I_2)_{a}\vert \leq C_\eta e^{5\eta t}}$. \\
\\
The term $(I_2)_b$ is estimated following the same computations done for the term $(I_2)_a$ above and thus we have the estimate (\ref{eq94}). \\
\\
Finally, with the estimate (\ref{eq94}) we get back to the estimate (\ref{eq95}) and we write 
\begin{equation}\label{eq96}
\vert I_2 \vert \leq C_\eta \frac{e^{5 \eta t }}{\vert x \vert^3},
\end{equation}  and thus,  by the estimates (\ref{eq97}) and (\ref{eq96}) at hand, we get back to the estimate (\ref{eq98}) and we can write the desired inequality:  $\ds{\vert \partial_x K_{t,x}\vert \leq C_\eta \frac{e^{5 \eta t }}{\vert x \vert^3}}$. 
\item[2)]  We write 
\begin{equation}
\vert \partial_x K_\eta(t,x)\vert \leq  \int_{\R} \vert (2\pi i \xi) e^{2\pi i x \xi} \widehat{K_{\eta}}(t,\xi) \vert d \xi \leq \Vert (1+ \vert \xi \vert) \widehat{K_{\eta}}(t,\cdot)\Vert_{L^1},
\end{equation}  and by the estimate (\ref{Estim-m}) (with $m=1$) we have 
\begin{equation}
\Vert (1+ \vert \xi \vert) \widehat{K_{\eta}}(t,\cdot)\Vert_{L^1} \leq C_\eta \left(e^{2 \eta t} + \frac{1}{t^{\frac{1}{3}}}+ \frac{1}{t^{\frac{2}{3}}} \right)= \frac{C_\eta }{t^{\frac{2}{3}}} \left( t^{\frac{2}{3}} e^{2 \eta t} + t^{\frac{1}{3}} +1\right) \leq \frac{C_\eta}{t^{\frac{2}{3}}} e^{5 \eta t}. 
\end{equation} Then we can write 
\begin{equation*}
\vert \partial_x K_\eta(t,x) \vert \leq \frac{C_\eta}{t^{\frac{2}{3}}} e^{5 \eta t}\leq \frac{C_\eta}{t^{\frac{2}{3}}} e^{6 \eta t}.
\end{equation*} Finally, by this estimate and the estimate given in point $1)$ above:  $\ds{\vert \partial_x K_{t,x}\vert \leq C_\eta \frac{e^{5 \eta t }}{\vert x \vert^3}}$, we obtain: $\ds{\vert \partial_x K_\eta(t,x)\vert} \leq \ds{C_\eta \frac{e^{6\eta t}}{t^{\frac{2}{3}}} \frac{1}{1+\vert x \vert^3}}$.  \finpv
\end{enumerate}	
%%%%%%%%%%%%%%%%%%%%%%%%%%%%%%%%%%%%%%%%%%%

\end{document}